\documentclass[a4paper,11pt]{amsart} %

\usepackage{geometry}
\usepackage{hyphenat} %
\usepackage[utf8]{inputenc}
\usepackage[T1]{fontenc}
\usepackage[all,cmtip,2cell]{xy}
\UseTwocells
\xyoption{2cell}
\usepackage{amsfonts}
\usepackage{amsthm}
\usepackage{amsmath}
\usepackage{amssymb}
\usepackage{bbm} %
\usepackage{enumerate}
\usepackage{verbatim} %
\usepackage{mathtools}
\usepackage[usenames,dvipsnames]{xcolor}
\usepackage{url}
\usepackage{chngcntr} %
\usepackage{pdflscape} %

\usepackage[titletoc]{appendix}

\usepackage{tikz-cd}

\usepackage{etoolbox}
\makeatletter
\appto{\appendices}{\def\Hy@chapapp{Appendix}}
\makeatother

\usepackage[stretch=10]{microtype}

\usepackage[pdftex]{hyperref} %
\usepackage{bookmark}

\renewcommand{\epsilon}{\varepsilon}

\theoremstyle{definition} \newtheorem{defn}[equation]{Definition}
\theoremstyle{remark} \newtheorem{notation}[equation]{Notation}
\theoremstyle{plain} \newtheorem{teo}[equation]{Theorem}
\theoremstyle{plain} \newtheorem{lema}[equation]{Lemma}
\theoremstyle{plain} \newtheorem{prop}[equation]{Proposition}
\theoremstyle{plain} \newtheorem{corolario}[equation]{Corollary}
\theoremstyle{remark} \newtheorem{obs}[equation]{Remark}
\theoremstyle{remark} 
\theoremstyle{remark} 
\theoremstyle{definition} 
\theoremstyle{remark} \newtheorem{ej}[equation]{Example}
\theoremstyle{remark} 
\theoremstyle{plain}

\DeclareUnicodeCharacter{00A0}{ } %

\DeclareMathAlphabet{\mathbbe}{U}{bbold}{m}{n}

\newcommand{\bq}{\begin{quote} \it}
\newcommand{\eq}{\it\end{quote}}

\newcommand{\paragrafo}[1]{\textbf{#1}\ \ \ }

\newcommand{\F}{\mathbb{F}}

\renewcommand{\L}{\mathcal{L}}

\renewcommand{\1}{\ensuremath{\mathbbm{1}}}
\newcommand{\V}{\mathcal{V}}

\newcommand{\W}{\mathcal{W}}

\newcommand{\Q}{\mathbb{Q}}
\renewcommand{\S}{\mathbb{S}}

\newcommand{\Z}{\mathbb{Z}}
\newcommand{\N}{\mathbb{N}}

\newcommand{\be}{\begin{enumerate}}
\newcommand{\ee}{\end{enumerate}}
\newcommand{\bi}{\begin{itemize}}
\newcommand{\ei}{\end{itemize}}
\newcommand{\bprop}{\begin{prop}}
\newcommand{\eprop}{\end{prop}}
\newcommand{\bej}{\begin{ej}}
\newcommand{\eej}{\end{ej}}
\newcommand{\bdefn}{\begin{defn}}
\newcommand{\edefn}{\end{defn}}
\newcommand{\bprf}{\begin{proof}}
\newcommand{\eprf}{\end{proof}}
\newcommand{\bobs}{\begin{obs}}
\newcommand{\eobs}{\end{obs}}
\newcommand{\bcor}{\begin{corolario}}
\newcommand{\ecor}{\end{corolario}}
\newcommand{\bteo}{\begin{teo}}
\newcommand{\eteo}{\end{teo}}

\newcommand{\id}{\mathrm{id}}

\newcommand{\Mon}{\ensuremath{\mathbf{Mon}}}
\newcommand{\CMon}{\ensuremath{\mathbf{CMon}}}

\newcommand{\Mod}{\mbox{-}\mathbf{Mod}}

\newcommand{\CAlg}{\mbox{-}\mathbf{CAlg}}

\newcommand{\RMod}{\ensuremath{R}\mbox{-}\mathbf{Mod}}
\newcommand{\AMod}{\ensuremath{A}\mbox{-}\mathbf{Mod}}

\newcommand{\SCAlg}{\ensuremath{\mathbb{S}}\mbox{-}\ensuremath{\mathbf{CAlg}}}
\newcommand{\RCAlg}{\ensuremath{R}\mbox{-}\ensuremath{\mathbf{CAlg}}}

\newcommand{\ACAlg}{\ensuremath{A}\mbox{-}\ensuremath{\mathbf{CAlg}}}

\newcommand{\SMod}{\ensuremath{\mathbb{S}}\mbox{-}\ensuremath{\textbf{Mod}}}

\newcommand{\Top}{\ensuremath{\mathbf{Top}}}

\newcommand{\Tor}{\mathrm{Tor}}

\newcommand{\cy}{\mathrm{cy}}
\newcommand{\sip}{\Sigma^\infty_+}

\renewcommand{\top}{\mathrm{top}}

\newcommand{\bigslant}[2]{{\raisebox{.2em}{$#1$}\left/\raisebox{-.2em}{$#2$}\right.}}

\newcommand{\paralelas}[2]{\ar@<.2pc>[r]^-{#1} \ar@<-.2pc>[r]_-{#2}}
\newcommand{\paralelass}[2]{\ar@<.2pc>[rr]^-{#1} \ar@<-.2pc>[rr]_-{#2}}

 \numberwithin{equation}{section}

\newcommand{\K}{\mathbf{K(\Z,2)}}
\renewcommand{\L}{\mathcal{L}}

\newcommand{\sma}{\wedge}

\hypersetup{pdfauthor={Bruno Stonek},colorlinks=true, citecolor=blue, urlcolor=blue, linkcolor=black,breaklinks=true,hypertexnames=false}

\author{Bruno Stonek}
\address{Max Planck Institute for Mathematics, Vivatsgasse 7, 53111 Bonn, Germany.}
\email{stonek@mpim-bonn.mpg.de} 

\begin{document}

\title{Higher topological Hochschild Homology of periodic complex K-theory}

\def\C{\mathbb{C}}

\begin{abstract} We describe the topological Hochschild homology of the periodic complex $K$-theory spectrum, $THH(KU)$, as a commutative $KU$-algebra: it is equivalent to $KU[K(\Z,3)]$ and to $F(\Sigma KU_\Q)$, where $F$ is the free commutative $KU$-algebra functor on a $KU$-module. Moreover, $F(\Sigma KU_\Q)\simeq KU \vee \Sigma KU_\Q$, a square-zero extension. In order to prove these results, we first establish that topological Hochschild homology commutes, as an algebra, with localization at an element. 

Then, we prove that $THH^n(KU)$, the $n$-fold iteration of $THH(KU)$, i.e. $T^n\otimes KU$, is equivalent to $KU[G]$ where $G$ is a certain product of integral Eilenberg-Mac Lane spaces, and to a free commutative $KU$-algebra on a rational $KU$-module.
We prove that $S^n \otimes KU$ is equivalent to $KU[K(\Z,n+2)]$ and to $F(\Sigma^n KU_\Q)$. We describe the topological André-Quillen homology of $KU$ as $KU_\Q$.
\end{abstract}

\maketitle

\section{Introduction}

Topological Hochschild homology ($THH$) of structured ring spectra was introduced by Bökstedt \cite{bokstedt} and Breen \cite{breen}; for an introduction to the subject, see \cite[Chapter 4]{dgm}, \cite[Chapter IX]{ekmm} and \cite{shipley-thh}. It is the generalization to structured ring spectra of classical Hochschild homology ($HH$) of rings.

It was realized in \cite{mc-schw-vo} that the $THH$ of a commutative ring spectrum $R$ can be expressed as $S^1\otimes R$, where $\otimes$ denotes the tensor of the category of commutative ring spectra over unbased spaces. Tensors with other spaces also give interesting information \cite{bcd}, \cite{cdd}, and we refer to them as giving ``higher $THH$'' of $R$. For example, tori $T^n$ give $n$-fold iterated $THH$. %
Spheres $S^n$ give a topological version of Pirashvili's higher order Hochschild homology of commutative rings \cite{pirashvili}. Complete calculations of these invariants for a given $R$ are scarce: see for example \cite{schlichtkrull-higher} for the case of  spectra, \cite{veen} and \cite{blprz} for partial computations for the Eilenberg-Mac Lane ring spectrum $H\F_p$ of the field with $p$ elements  and other related ring spectra, and \cite{dlr} for Eilenberg-Mac Lane spectra of some rings of integers.

In this paper, we present complete descriptions of the commutative $KU$-algebras $T^n\otimes KU$ and $S^n\otimes KU$ for $n\geq 1$, where $KU$ is the ring spectrum of periodic complex topological $K$-theory. Close results to the $n=1$ case were known, as we explain below, but only additively. The classical André-Quillen homology of commutative rings also has a topological analogue, denoted $TAQ$ \cite{basterra}: we determine the $KU$-module $TAQ(KU)$.

Our computations showcase some interesting phenomena. The formulas for the (higher) topological Hochschild homology and topological André-Quillen homology of $KU$ which we obtain are the ones we would get if $KU$ was somehow a Thom spectrum (which it isn't), 
see Remarks \ref{thh-thom}, \ref{x-thom} and \ref{jurs}.\ref{taq-ku-thom}. Also, our results show that the conclusion of McCarthy-Minasian's adaptation of the Hochschild-Kostant-Rosenberg theorem \cite{mccarthy-minasian}, which applies only to connective ring spectra, holds for $KU$ (Remark \ref{jurs}.\ref{hokr}).

Finally, another remarkable phenomenon highlighted by our computations involves invariance under stable equivalences of spaces: let $R$ be a commutative ring spectrum and $X$ and $Y$ be spaces such that $\Sigma X \simeq \Sigma Y$. One may ask the question of whether $X\otimes R\simeq Y\otimes R$, i.e. of whether $-\otimes R$ is a \emph{stable invariant}. This turns out not to be true in general \cite{dundas-tenti}, but the computations of \cite{veen} show that, in a certain range relating $n$ and $p$, $X\otimes H\F_p \simeq Y\otimes H\F_p$ when $X=T^n$ and $Y=\bigvee\limits_{i=1}^n (S^i)^{\vee {n \choose i}}$. Our results show that the same is true for $R=KU$, see Remark \ref{stabeq}. It would be interesting to know whether $-\otimes KU$ is a stable invariant, and to try to characterize the commutative ring spectra $R$ such that $-\otimes R$ is a stable invariant.\\

\paragrafo{Summary of results} We work in the context of $\S$-modules and commutative $\S$-algebras from \cite{ekmm}. We use an adaptation of the model for $KU$ given by Snaith \cite{snaith79}, \cite{snaith81}, namely $\sip K(\Z,2)[x^{-1}]$, to this context. If we are to use this model to compute $THH(KU)$, we first need to prove that $THH$ commutes with localizations: this is done in \textbf{Corollary \ref{cor-thhloc}}, taking care of the multiplicative structure.

Our first expression for $THH(KU)$ as a commutative $KU$-algebra is obtained in \textbf{Theorem \ref{thhku1}}:
\[THH(KU)\simeq KU[K(\Z,3)],\]
where the underlying $KU$-module of $KU[K(\Z,3)]$ is $KU \wedge K(\Z,3)_+$. The second one is given in \textbf{Theorem \ref{thhku2}}: there are weak equivalences of commutative $KU$-algebras
\[\xymatrix{ KU \vee \Sigma KU_\Q  & \ar[l]_-\sim F(\Sigma KU_\Q) \ar[r]^-\sim & THH(KU).}\]
Here $F$ denotes the free commutative $KU$-algebra on a $KU$-module functor, and the $KU$-algebra structure on $KU\vee \Sigma KU_\Q$ is that of a square-zero extension. 

Note that, previously, McClure and Staffeldt \cite[8.1]{mc-st} established that $THH(L)\simeq L \vee \Sigma L_\Q$ as spectra, where $L$ is the $p$-adic completion of the Adams summand of $KU$ for a given odd prime $p$. In \cite[7.9]{thhko}, the authors show that $THH(KO)\simeq KO\vee \Sigma KO_\Q$ as $KO$-modules; here $KO$ is the periodic real complex $K$-theory ring spectrum. Note that these results (and others closely related, see Remark \ref{antecedentes}) are about the additive structure and do not involve the multiplicative structure, which we take into consideration. Finally, note that a lot of effort was devoted to describing $THH(ku)$ \cite{ausoni-thhku}, where $ku$ denotes the connective complex $K$-theory spectrum: that case is markedly harder.\\

We consider the iterated $THH$ of $KU$. The first expression we gave above for $THH(KU)$ directly generalizes: one replaces $K(\Z,3)$ by a suitable product of integral Eilenberg-Mac Lane spaces. See \textbf{Theorem \ref{thhnku1cor}}: there is a zig-zag of weak equivalences of commutative $KU$-algebras
\[THH^n(KU)\simeq KU\left[\prod\limits_{i=1}^n K(\Z,i+2)^{\times {n \choose i}}\right].\]
The second expression for $THH(KU)$ also generalizes: this is \textbf{Theorem \ref{thhnkufree}}, where we get a zig-zag of weak equivalences of commutative $KU$-algebras
\[F\left(\bigvee\limits_{i=1}^n (S^i)^{\vee {n\choose i}} \wedge KU_\Q\right) \simeq T^n\otimes KU.\]
The expression $KU \vee \Sigma KU_\Q$ for $THH(KU)$ also generalizes to $THH^n(KU)$. %
In this case, the augmentation ideal $\overline{THH}^n(KU)$ is still rational, but it has a non-trivial non-unital commutative $KU$-algebra structure. We describe the non-unital commutative $\Q[t^{\pm 1}]$-algebra $\overline{THH}^n_*(KU)$ as iterated Hochschild homology. See \textbf{Theorem \ref{high-hoch}}.%
\\

We then shift our attention to $X\otimes KU$, where $X$ is a based CW-complex which is a reduced suspension, e.g. a sphere $S^n$. In this case, the first description for $THH(KU)$ generalizes as a zig-zag of weak equivalences of commutative $KU$-algebras:
\[S^n\otimes KU \simeq KU[K(\Z,n+2)].\]
This is \textbf{Theorem \ref{holo}}. The second description for $THH(KU)$ generalizes as %
\[F(S^n \wedge KU_\Q) \simeq S^n\otimes KU.\] 
This is a particular case of \textbf{Theorem \ref{xku}}. Finally, we establish that
\[TAQ(KU)\simeq KU_\Q\]
as $KU$-modules, where $TAQ(KU)$ denotes the topological André-Quillen $KU$-module of $KU$. This is \textbf{Corollary \ref{cortaq}}.\\

\paragrafo{Remark on trace methods} One reason for the importance of $THH$ is its relation to algebraic $K$-theory. If $R$ is a (discrete) ring, then the trace map $K(R)\to HH(R)$ factors through the topological Hochschild homology of the Eilenberg-Mac Lane ring spectrum of $R$. Moreover, the trace map $K(A)\to THH(A)$ exists for any ring spectrum $A$. Out of topological Hochschild homology one can build topological cyclic homology, $TC(A)$, and the trace map further factors through it. The spectrum $TC(A)$ is closely related to algebraic $K$-theory: see \cite{dgm}. We might thus see $THH$ as a more easily approachable stepping stone on the way to the more fundamental algebraic $K$-theory.

It is important to note, however, that $THH(KU)$ is unlikely to be of assistance in the determination of $K(KU)$ via the methods we pointed out in the previous paragraph. First of all, note that one of the most useful theorems for computing algebraic $K$-theory via trace maps, namely, the theorem of Dundas-Goodwillie-McCarthy \cite[7.0.0.2]{dgm} only applies to connective ring spectra. Therefore, one may wish to get to $K(KU)$ by noting that $KU$ is the localization of $ku$, and by applying trace methods for $ku$. Indeed, in \cite{bm-localization-08}, Blumberg and Mandell prove that $K(KU)$ sits in a localization cofiber sequence $K(\Z)\to K(ku)\to K(KU)\to \Sigma K(\Z)$. In \cite{bm-localization-14}, they establish an analogous cofiber sequence for $THH$, but the term involving $KU$ is \emph{not} $THH(KU)$ but an appropriate modification of it which receives a trace map from $K(KU)$; Ausoni shows in \cite[8.3]{ausoni-2010} how to compute the $V(1)$-homotopy ($p$ odd) of $K(KU)$ using this approach, and in \cite[3.6]{ausoni-rognes-rational} him and Rognes determine $K(KU)$ rationally. Note that these computations do not involve $THH(KU)$.

On the other hand, an interesting question to ask is if there are any elements in $K(KU)$ which survive to $THH(KU)$ via the trace. 
We know that that $V(1)_*THH(KU)=0$ and $V(0)_*THH(KU)\cong V(0)_*KU$, but rationally, the trace $K(KU)\to THH(KU)$ is non-zero (see \cite[Paragraph 5.3]{ausoni-rognes-rational}). Therefore, the trace $K(KU)\to THH(KU)$ might well be useful in studying the integral homotopy type of $K(KU)$. 
More generally, it would be interesting to detect elements in $THH^n(KU)$ that survive from the $n$-fold iterated algebraic $K$-theory of $KU$ via the iterated trace. See \cite{cdd}: they propose $THH^n$ as ``a computationally tractable cousin of $n$-fold iterated algebraic $K$-theory''.\\

\paragrafo{Comment on cofibrancy} At the heart of the computation of $THH(KU)$ lies the isomorphism $THH(\S[G])\cong \S[B^\cy G]$ where $G$ is a topological commutative monoid. The core of this result was already known; we prove it in Proposition \ref{conmutarTHH}, taking care of the multiplicative structures. This is a point-set result. However, the procedure of localization of a commutative $\S$-algebra at a homotopy element takes a \emph{cofibrant} commutative $\S$-algebra as input, so if we are to exploit the isomorphism just stated for the computation of $THH(KU)$ via Snaith's theorem, we first need to prove that $THH$ preserves the weak equivalence to $\S[G]$ from a cofibrant commutative $\S$-algebra replacement of it. This is obtained in Section \ref{sect:cofibrancy}. Assuming $G$ is a CW-complex with unit a 0-cell, the key property that $\S[G]$ satisfies is that it is \emph{flat}, i.e. smashing an $\S$-module with it preserves weak equivalences. This, along with the fact that the simplicial cyclic bar construction $B^\cy_\bullet\S[G]$ is a \emph{proper} simplicial $\S$-module, proves to be enough.\\

\paragrafo{Outline of the paper}  In Section \ref{sect-cof}, we review some model categorical aspects of \cite{ekmm}, particularly those pertaining to commutative $\S$-algebras. In Section \ref{inversion}, we prove some elementary properties of localization of a commutative $\S$-algebra at an element. In Section \ref{sect-thh}, we review some needed aspects of topological Hochschild homology, and we prove that $THH$ commutes with localization at an element. Section \ref{sect:thhku} contains the results pertaining to $THH(KU)$, %
and in Sections \ref{sect-iterated} and \ref{section:snku} we prove our results about $T^n\otimes KU$ and $S^n\otimes KU$. Finally, in Section \ref{section:taq}, we determine the topological André-Quillen homology of $KU$.\\

\paragrafo{Conventions} By \emph{space} we will mean ``compactly generated weakly Hausdorff topological space'', and we will denote the cartesian closed category they form by $\Top$, which we endow with the Quillen model structure. The corresponding model category of based spaces will be denoted by $\Top_*$. 
We will work with the categories of \cite{ekmm}: our main objects are $\S$-modules, commutative $\S$-algebras $R$, $R$-modules and commutative $R$-algebras $A$.\\

\paragrafo{Acknowledgments} I would like to thank Christian Ausoni, my PhD supervisor, for suggesting this project, sharing his ideas and his support; Geoffroy Horel for his very useful suggestions; Eva Höning for our many engaging and fruitful discussions; Christian Schlichtkrull for his careful reading and his corrections; and Bj\o rn Dundas for so warmly and selflessly sharing so much of his time, ideas and insights at the University of Bergen. I would also like to thank Tobias Barthel and Magdalena Zielenkiewicz for their assistance with later revisions.

Parts of the content of this article are part of the author's PhD dissertation at Université Paris 13. Research partially supported by the \emph{ANR-16-CE40-0003 project ChroK} and the \emph{Fondation Sciences Mathématiques de Paris} (FSMP). The author would like to thank the Max Planck Institute for Mathematics at Bonn for their hospitality. 

\section{Model structures} \label{sect-cof}

We will freely use the language of (enriched, monoidal) model categories as expounded in e.g. \cite{mayponto}.

The category $\SMod$ of $\S$-modules has a $\Top_*$-enriched symmetric monoidal cofibrantly generated model structure \cite[VII.4]{ekmm}. %
A commutative $\S$-algebra is, by definition, a commutative monoid in $\SMod$. The category they form, $\SCAlg$, can also be described as the category of $\mathbb P$-algebras where $\mathbb P$ is the commutative monoid monad. The forgetful functor $U:\SCAlg\to \SMod$ creates a model structure on $\SCAlg$\footnote{A functor $U:\mathcal C\to \mathcal M$ \emph{creates a model structure} on $\mathcal C$ if $\mathcal M$ is a model category and $\mathcal C$ is a model category such that $f$ is a fibration (resp. weak equivalence) in $\mathcal C$ if and only if $Uf$ is a fibration (resp. weak equivalence) in $\mathcal M$. We say that $U$ \emph{strongly creates} the model structure of $\mathcal C$ if, in addition, $f$ is a cofibration in $\mathcal C$ if and only if $Uf$ is a cofibration in $\mathcal M$. We are following the nomenclature of \cite[15.3.5]{mayponto}.}. In particular, there is a Quillen adjunction $\xymatrix{\SMod\ar@<.2pc>[r]^-F  &   \ar@<.2pc>[l]^-U \SCAlg}$. The category $\SCAlg$ has a $\Top$-enriched symmetric monoidal cofibrantly generated model structure. In both $\SMod$ and $\SCAlg$, the monoidal product is the smash product $\wedge$ and the monoidal unit is the sphere spectrum $\S$.\footnote{Note that the pushout-product axiom is satisfied in $\SCAlg$: indeed, the smash product of commutative $\S$-algebras $R$ and $T$ is their coproduct, i.e. the pushout of $T \leftarrow \S \to R$ \cite[II.3.7]{ekmm}. Therefore, a ``pushouts commute with pushouts'' argument proves that the pushout-product map is an isomorphism. This proves, more generally, that any model category with the cocartesian monoidal structure (i.e. the monoidal product is the coproduct and the unit is the initial object) satisfies the pushout-product axiom and hence is a monoidal model category.}

Let $R\in \SCAlg$, and consider the category of $R$-modules, $\RMod$. The forgetful functor $\RMod\to \SMod$ creates a model structure on $\RMod$, and $\RMod$ acquires a $\Top_*$-enriched symmetric monoidal cofibrantly generated model category structure. The forgetful functor $U:\RCAlg\to \RMod$ creates a model structure on $\RCAlg$, and thus $\RCAlg$ has a $\Top$-enriched symmetric monoidal cofibrantly generated model category structure. In both $\RCAlg$ and $\RMod$, the monoidal product is the smash product relative to $R$, $\wedge_R$, and the monoidal unit is $R$. In all these model categories, all objects are fibrant.

Cofibrancy is more delicate. The sphere $\S$-module $\S$ is not cofibrant as an $\S$-module, but it is cofibrant as a commutative $\S$-algebra.  More generally, the underlying $R$-module of a cofibrant commutative $R$-algebra is generally not cofibrant as an $R$-module.

Let $R$ be a commutative $\S$-algebra. We record the following useful properties:

\be
\item \label{ujin} \label{moguil} If $M$ is a cofibrant $R$-module, then $M\wedge_R -$ preserves all weak equivalences of $R$-modules, so if $X$ is any $R$-module, then $M\wedge_R X$ represents the derived smash product \cite[III.3.8]{ekmm}. Note that $X\wedge_R-$ preserves weak equivalences between cofibrant $R$-modules. Indeed: let $f:M\to N$ be such a weak equivalence. Let $\gamma_X:\Gamma X\to X$ be a cofibrant replacement of $X$. We have a commutative diagram
\begin{equation}\label{moguil-diag}\xymatrix@C+1pc{X\wedge_R M \ar[r]^-{\id \wedge f} & X\wedge_R N  \\ \Gamma X \wedge_R M \ar[r]^-\sim_-{\id \wedge f} \ar[u]_-\sim^-{\gamma_X \wedge \id} & \Gamma X \wedge_R N \ar[u]^-\sim_-{\gamma_X \wedge \id}}\end{equation}
where the two vertical maps and the bottom horizontal map are weak equivalences by the result just quoted, so the top vertical map is a weak equivalence, too.
\item \label{es} Suppose $R$ is cofibrant as a commutative $\S$-algebra. %
Let $A$ and $B$ be cofibrant commutative $R$-algebras. Let $\gamma_A:\Gamma A \to A$ and $\gamma_B: \Gamma B \to B$ be cofibrant replacements of $A$ and $B$ in the category of $R$-modules. Then \[\gamma_A \wedge \gamma_B: \Gamma A \wedge_R \Gamma B \to A\wedge_R B\] is a weak equivalence of $R$-modules \cite[VII.6.4, 6.5, 6.7]{ekmm}. This tells us that $A\wedge_R B$ computes the derived smash product of $A$ and $B$ as $R$-modules. As a consequence of this and of (\ref{ujin}), by the 2-out-of-3 property we deduce that \[\gamma_A \wedge \id_B: \Gamma A \wedge_R B \to A \wedge_R B\] is a weak equivalence, since $\gamma_A \wedge \gamma_B = (\gamma_A \wedge \id_B) \circ (\id_{\Gamma A} \wedge \gamma_B)$. Similarly, $\id_A \wedge \gamma_B$ is also a weak equivalence.
\item As in any model category, the coproduct of cofibrant objects is cofibrant. Hence, if $A$ and $B$ are cofibrant commutative $R$-algebras, then $A\wedge_R B$ is a cofibrant commutative $R$-algebra \cite[VII.6.8]{ekmm}.
\item Let $\S\to A\to B$ be cofibrations of commutative $\S$-algebras. Then the functor $B\wedge_A-:A\CAlg\to B\CAlg$ preserves weak equivalences between commutative $A$-algebras which are cofibrant as commutative $\S$-algebras \cite[VII.7.4]{ekmm}.
\item  \label{cinco} The category $\RCAlg$ can also be described as the category of objects of $\SCAlg$ under $R$. As such, the forgetful functor $\RCAlg \to \SCAlg$ strongly creates a model structure on $\RCAlg$ \cite[Theorem 15.3.6]{mayponto}. This model structure coincides with the one described above \cite[Remark 2.4.1]{eva-thesis}. In conclusion, a map $f:A\to B$ is a cofibration in $\RCAlg$ if and only if it is a cofibration in $\SCAlg$. In particular, if $R$ is a cofibrant commutative $\S$-algebra and $A$ is a cofibrant commutative $R$-algebra, then $A$ is cofibrant as a commutative $\S$-algebra. 
\ee

Note: in \cite{ekmm} they call \emph{q-cofibration} what we call a cofibration. We will have no use for what they call a cofibration. We use the term ``homotopy category'' in the model categorical sense \cite[14.4.1]{mayponto}: in \cite{ekmm} these were called \emph{derived categories}. 

\bobs\label{funct-fact} The following remarks will be used below. In a cofibrantly generated model category constructed via Quillen's small object argument, the factorization of an arrow as a cofibration followed by an acyclic fibration is \emph{functorial} \cite[12.2.2]{riehl}, a concept carefully defined e.g. in \cite[12.1.1]{riehl}, \cite[14.1.10]{mayponto}. As an example, the category $\SCAlg$ admits such a functorial factorization. We can factor the unit maps in order to obtain a functorial cofibrant replacement functor $Q$.

Let $R$ be a commutative $\S$-algebra. If $A$ is a commutative $R$-algebra, then we may factor the unit map $\eta:R\to A$ with the functorial factorization in $\SCAlg$ described above. Denote by $Q_RA$ the object appearing in the factorization, i.e. $\eta$ is factored as $R\to Q_RA\to A$. From item (\ref{cinco}) above we obtain that the first arrow is a cofibration in $\RCAlg$ and the second arrow is an acyclic fibration in $\RCAlg$. The functoriality of the factorization proves that this defines a functor $Q_R$ of cofibrant replacement in the category $\RCAlg$.
\eobs

We end this section on model structures with the following general lemma which will be used in the proof of Theorem \ref{xku}. We thank Eva Höning for explaining this lemma to us.

\begin{lema} \label{hompus} Let \begin{equation}\label{tef}\xymatrix@C+1pc{B \ar[d]_-u^-\sim & A \ar[l]_-f \ar[d]_-\sim^-v \ar[r]^-g & C \ar[d]^-w_-\sim \\ B' & A' \ar[l]^-{f'} \ar[r]_ -{g'} & C'}\end{equation}
be a diagram in a model category where the vertical arrows are weak equivalences, all objects are cofibrant and one map in each horizontal line is a cofibration. Suppose both squares are homotopy commutative. Then there is a natural zig-zag of weak equivalences between the pushouts of both horizontal lines.%
\bprf  %
Let $(\textup{Cyl}(A),i_0,i_1)$ denote a cylinder object for $A$. Let $H:\textup{Cyl}(A)\to B'$ denote a homotopy from $uf$ to $f'v$, and $G:\textup{Cyl}(A)\to C'$ denote a homotopy from $wg$ to $g'v$. We have the following (strictly) commutative diagram:
\[\xymatrix{
B \ar[d]_-u^-\sim & A \ar[l]_-f \ar[r]^-g \ar[d]^-\id & C \ar[d]^-w_-\sim \\
B' \ar[d]_-\id & A \ar[l]_-{uf} \ar[r]^-{wg} \ar[d]^-{i_0}_-\sim & C' \ar[d]^-\id \\
B'  & \textup{Cyl}(A) \ar[l]_-H \ar[r]^-G & C' \\
B' \ar[d]_-\id \ar[u]^-\id & A \ar[d]^-v_-\sim \ar[l]_-{f'v} \ar[r]^-{g'v} \ar[u]_-{i_1}^-\sim & C' \ar[d]^-\id \ar[u]_-\id \\
B' & A' \ar[l]^-{f'} \ar[r]_-{g'} & C' }\]
By a repeated application of the homotopy invariance of homotopy pushouts \cite[Dual of 13.3.4]{hirschhorn}, we obtain a zig-zag of weak equivalences between the homotopy pushout of $(f,g)$ and the homotopy pushout of $(f',g')$. But these homotopy pushouts are computed by the (categorical) pushouts, since in (\ref{tef}) all objects are cofibrant and one map in each line is a cofibration \cite[Dual of 13.3.8]{hirschhorn}.
\eprf
\end{lema}

\section{Inversion of an element}\label{inversion} %

In this section, we recall the procedure of inverting a homotopy element in a commutative $\S$-algebra following \cite{ekmm} and prove some properties which will be needed below.

\bteo \label{rect} \cite[VIII.2.2, VIII.4.2]{ekmm} Let $R$ be a cofibrant commutative $\S$-algebra and $x\in \pi_*R$. There exists a cofibrant %
commutative $R$-algebra $R[x^{-1}]$ with unit $j:R\to R[x^{-1}]$ satisfying that $\pi_*(R[x^{-1}])=\pi_*(R)[x^{-1}]$, and if $f:R\to T$ is a map in $\SCAlg$ such that $(\pi_*f)(x)\in \pi_*T$ is invertible, then there exists a map $\tilde f:R[x^{-1}]\to T$ in $\SCAlg$ making the following diagram commute:
\[\xymatrix{R\ar[r]^-f \ar[d]_-j & T \\ R[x^{-1}] \ar@{.>}[ru]_{\tilde{f}} }.\]
The map $\tilde f$ is unique up to homotopy of commutative $\S$-algebras. 
Moreover, if the morphism $\pi_*(R)[x^{-1}]\to \pi_*T$ coming from the universal property for localizations of commutative $\pi_*(R)$-algebras is an isomorphism, then $\tilde f$ is a weak equivalence.
\eteo

The previous theorem is valid, \emph{mutatis mutandis}, if $\S$ is replaced by some cofibrant commutative $\S$-algebra.

\begin{lema} \label{idem} The multiplication map $\mu:R[x^{-1}]\wedge_R R[x^{-1}]\to R[x^{-1}]$ is a weak equivalence of commutative $R[x^{-1}]$-algebras.
\bprf The Tor spectral sequence \cite[IV.4.1]{ekmm} here takes the form %
\[E^2_{*,*}=\Tor_{*,*}^{\pi_*R}(\pi_* R[x^{-1}], \pi_* R[x^{-1}]) \Rightarrow \pi_*(R[x^{-1}]\wedge_R R[x^{-1}]). \] %
Since the localization morphism $\pi_*R\to \pi_*R[x^{-1}]$ is flat, the spectral sequence is concentrated in the 0-th column and thus the edge homomorphism
\begin{equation}\label{edge} \nabla: \pi_*R[x^{-1}]\otimes_{\pi_*R}\pi_*R[x^{-1}] \to \pi_*(R[x^{-1}]\wedge_R R[x^{-1}])\end{equation}
is an isomorphism. Since $\wedge_R$ is the coproduct in the category of commutative $R$-algebras, we can consider  the canonical maps $i_1,i_2:R[x^{-1}]\to R[x^{-1}]\wedge_R R[x^{-1}]$. The edge homomorphism $\nabla$ coincides with the map $(\pi_*i_1, \pi_*i_2)$ defined via the universal property of the coproduct of commutative $\pi_*R$-algebras. %
We have the following commutative diagram of commutative $\pi_*R$-algebras:
\[\xymatrix{
\pi_*R[x^{-1}] \ar@/_2pc/[ddr]_-\id \ar[r]^-{\iota_1} \ar[rd]_-{\pi_*i_1}
& \pi_*R[x^{-1}]\otimes_{\pi_*R} \pi_*R[x^{-1}] \ar[d]^-\nabla
& \pi_*R[x^{-1}] \ar[l]_-{\iota_2} \ar[ld]^-{\pi_*i_2} \ar@/^2pc/[ddl]^-\id
\\
& \pi_*(R[x^{-1}] \wedge_R R[x^{-1}]) \ar[d]^-{\pi_*\mu}
\\
& \pi_*R[x^{-1}] }\] 
where $\iota_1,\iota_2$ are the canonical inclusions into a coproduct of commutative $\pi_*R$-algebras. Again, by the universal property of the coproduct of commutative $\pi_*R$-algebras, there is a unique arrow $\pi_*R[x^{-1}]\otimes_{\pi_*R} \pi_*R[x^{-1}] \to \pi_*R[x^{-1}]$ making the outer diagram commute. One such arrow is the canonical isomorphism that one has for any such algebraic localization, i.e. $h:S^{-1}A\otimes_A S^{-1}A \stackrel{\cong}{\to} S^{-1}A$ for any commutative ring $A$ and multiplicative subset $S\subset A$. Another such arrow is $\pi_*\mu \circ \nabla$. Therefore, $h=\pi_*\mu \circ \nabla$. Since $\nabla$ and $h$ are isomorphisms, so is $\pi_*\mu$.
\eprf
\end{lema}

If $f:R\to T$ is a morphism between cofibrant commutative $\S$-algebras and $x\in \pi_*R$, then Theorem \ref{rect} gives us a map of cofibrant commutative $\S$-algebras
\[\xymatrix{R\ar[r]^-f \ar[d]_-{j_R} & T \ar[d]^-{j_T} \\ R[x^{-1}] \ar@{.>}[r]_-{f[x^{-1}]} & T[(\pi_*f)(x)^{-1}]}\]
such that if $f$ is a weak equivalence, then $f[x^{-1}]$ is a weak equivalence. %
Note that $f[x^{-1}]$ turns $T[(\pi_*f)(x)^{-1}]$ into a commutative $R[x^{-1}]$-algebra. 

The previous square induces an arrow from the pushout $R[x^{-1}]\wedge_R T$ in $R\CAlg$. The following theorem tells us that it is a weak equivalence. Compare with \cite[V.1.15]{ekmm} which handles the case where $T$ is replaced by an $R$-module. 

\bprop[Base change for localization]\label{base-loc} Let $f:R\to T$ be a morphism of cofibrant commutative $\S$-algebras and $x\in \pi_*R$. The morphism of commutative $R$-algebras \begin{equation}\label{base-loceq}(f[x^{-1}],j_T): R[x^{-1}] \wedge_R T \to T[(\pi_*f)(x)^{-1}]\end{equation} %
is a weak equivalence.
\eprop
Note that (\ref{base-loceq}) is also a weak equivalence in $R[x^{-1}]\CAlg$ and in $T\CAlg$. Note as well that if $\epsilon:T\to R$ is a morphism of $\S$-algebras such that $\epsilon \circ f=\id_R$ so that $T$ becomes an augmented commutative $R$-algebra, then in (\ref{base-loceq}) both sides are naturally augmented over $R[x^{-1}]$ and the morphism commutes with the augmentations.

\bprf Denote the morphism $(f[x^{-1}],j_T)$ by $h$, for simplicity. Like in the proof of Lemma \ref{idem}, the Tor spectral sequence that computes the homotopy groups of $R[x^{-1}]\wedge_R T$ from those of $R[x^{-1}]$ and $T$ collapses, since $\pi_*R\to \pi_*R[x^{-1}]= (\pi_*R)[x^{-1}]$ is flat. %
Therefore, the map $\pi_*h$, fitting in a commutative diagram
\[\xymatrix{\pi_*(R[x^{-1}]\wedge_R T) \ar[r]^-{\pi_*h}  & (\pi_*T)[(\pi_*f)(x)^{-1}]\\ \ar[u]^-\cong (\pi_*R)[x^{-1}] \otimes_{\pi_*R} \pi_*T, \ar[ru]_-\cong  }\]
is an isomorphism, since the diagonal map is an isomorphism. Indeed, this is the map appearing in the analogous statement in commutative algebra of the theorem we are proving, applied to $\pi_*f: \pi_*R\to \pi_*T$. But this statement of commutative algebra is not hard to prove: it follows from the universal properties and the extension-restriction of scalars adjunction.
\eprf

\bprop \label{sepa} Let $R$ and $T$ be cofibrant commutative $\S$-algebras, $x\in \pi_nR$ and $y\in \pi_mT$. Denote by $x\wedge y$ the image of $x\otimes y$ under the morphism %
\[\xymatrix{\pi_*R \otimes_{\pi_*\S} \pi_*T \ar[r] & \pi_*(R\wedge T)}.\]
There is a weak equivalence of commutative $\S$-algebras 
 \[R[x^{-1}]\wedge T[y^{-1}] \to (R\wedge T)[(x\wedge y)^{-1}]\]
which is natural on $R$ and $T$.
\eprop %
Note that this is is also a map of commutative $R[x^{-1}]$ and $T[y^{-1}]$-algebras.
\bprf Let $i_1:R\to R\wedge T$, $i_2:T\to R\wedge T$ be the canonical maps into the coproduct. There exists a map $f$ making the following diagram commute.
\[\xymatrix{R \ar[d]_-{j_R} \ar[r]^{i_1} & R\wedge T \ar[r]^-{j_{R\wedge T}} & (R\wedge T)[(x\wedge y)^{-1}] \\ R[x^{-1}] \ar@{.>}[rru]_-f}\]
Indeed, applying $\pi_*$ to the horizontal composition, we get the map \[\pi_*(j_{R\wedge T}\circ i_1): \pi_*R \to \pi_*(R\wedge T)[(x\wedge y)^{-1}]\]which maps $x$ to $x\wedge 1$. This is an invertible element with inverse $(1\wedge y)(x\wedge y)^{-1}$, since the map $(\pi_*i_1,\pi_*i_2): \pi_*R \otimes_{\pi_*\S} \pi_* T \to \pi_*(R\wedge T)$ %
is multiplicative. Therefore, the property of Theorem \ref{rect} provides us with the arrow $f$ in $\SCAlg$. Similarly, we get a map $g:T[y^{-1}] \to (R\wedge T)[(x\wedge y)^{-1}]$. We assemble $f$ and $g$ into the coproduct map in $\SCAlg$
\[(f,g):R[x^{-1}]\wedge T[y^{-1}] \to (R\wedge T)[(x\wedge y)^{-1}].\]
Now recall from \cite[Section V.1]{ekmm} that $R[x^{-1}]$ is weakly equivalent, in $\RMod$, to the homotopy colimit of the tower
\begin{equation}\label{telescope}\xymatrix{R\ar[r]^-x & \Sigma^{-n} R \ar[r]^-x & \Sigma^{-2n} R\ar[r]^-x & \dots}.\end{equation}
The $T$-module $T[y^{-1}]$ is described similarly. The $R\wedge T$-module $(R\wedge T)[(x\wedge y)^{-1}]$ is weakly equivalent to the homotopy colimit of the tower
\[\xymatrix{R \wedge T\ar[r]^-{x\wedge y} & \Sigma^{-n-m} R\wedge T \ar[r]^-{x\wedge y} & \Sigma^{-2n-2m} R\wedge T\ar[r]^-{x\wedge y} & \dots}.\]
Smashing the homotopy colimit computing $R[x^{-1}]$ with the one computing $T[y^{-1}]$ we obtain the homotopy colimit computing $(R\wedge T)[(x\wedge y)^{-1}]$, since the diagonal map $\N\to \N\times \N$ is homotopy cofinal. The map $(f,g)$ is compatible with these identifications, hence it is a weak equivalence. %
\eprf

\begin{notation} Let $f:R\to T$ be a weak equivalence of cofibrant commutative $\S$-algebras. If $x\in \pi_*R$, we denote the algebra $T[(\pi_*f)(x)^{-1}]$ by $T[x^{-1}]$. Similarly, if $y\in \pi_*T$ we define $R[y^{-1}]$ as $T[(\pi_*f)^{-1}(y)^{-1}]$. In particular, if $R$ and $T$ are connected by a zig-zag of weak equivalences of cofibrant commutative $\S$-algebras, then a homotopy element $x\in \pi_*R$ defines $T[x^{-1}]$ and conversely. 
\end{notation}

We now turn to the inversion of an element in a non-cofibrant commutative $\S$-algebra $A$. Let $Q$ be a cofibrant replacement functor in the category of commutative $\S$-algebras, as obtained in Remark \ref{funct-fact}.%

\bdefn \label{def-nocof} Let $A\in \SCAlg$ and $x\in \pi_*A$. We define $A[x^{-1}]_h$ to be the cofibrant commutative $\S$-algebra $(QA)[x^{-1}]$.
\edefn

\bobs \label{lifto} If $\tilde A$ is a cofibrant commutative $\S$-algebra and $\tilde A\to A$ is a weak equivalence, then there is a weak equivalence of commutative $\S$-algebras $\tilde A\to QA$, and hence a weak equivalence $\tilde A[x^{-1}]\to QA[x^{-1}]$ of cofibrant commutative $\S$-algebras. Indeed, the existence of this weak equivalence follows from the lifting properties for the model category $\SCAlg$ and the 2-out-of-3 property for weak equivalences:
\[\xymatrix{\S \ar[r] \ar@{>->}[d] & QA \ar@{->>}[d]^-{\sim} \\ \tilde A \ar@{.>}[ru] \ar[r]_-\sim & A.}\]
In particular, the object $A[x^{-1}]_h$ in the homotopy category of commutative $\S$-algebras does not depend on the choice of $Q$, up to isomorphism. This explains the choice of the letter $h$ for the subscript.
\eobs

\section{Topological Hochschild homology} \label{sect-thh}

\subsection{Symmetric monoidal categories} 
Let $(\V,\otimes, \1)$ be a symmetric monoidal category. Denote by $\Mon(\V)$ and $\CMon(\V)$ the corresponding symmetric monoidal categories of monoids and commutative monoids in $\V$, respectively. Denote by $s\V$ the category of simplicial objects in $\V$; it is a symmetric monoidal category with levelwise monoidal product.

Suppose $\V$ is closed and cocomplete. Then if $A\in \CMon(\V)$, the category $A\Mod$ of $A$-modules gets a relative monoidal product $\otimes_A$ such that $(A\Mod, \otimes_A, A)$ is a symmetric monoidal category. One can thus speak of commutative $A$-algebras, which are the commutative monoids in $A\Mod$. We denote by $\ACAlg$ the category they form. %
An \emph{augmented} commutative $A$-algebra $B$ has an \emph{augmentation} map $B\to A$ which is a morphism of commutative monoids.

Let $F:\V\to \W$ be a strong symmetric monoidal functor between cocomplete closed symmetric monoidal categories, and suppose $F$ preserves colimits. Then $F$ induces a functor on commutative monoids, on modules over commutative monoids, and on commutative algebras. More specifically, there is an induced functor $F:\ACAlg \to F(A)\CAlg$. %

We will need the following
\begin{lema} \label{lemext} Let $F:\V\to \W$ be a functor as above, and let $A\in \CMon(\V)$. Then there is a natural isomorphism
\[\xymatrix@C+1pc{\CMon(\V) \ar[r]^-{A\otimes -} \ar[d]_-F & \ACAlg \ar[d]^-F \\ \CMon(\W) \ar[r]_-{F(A) \otimes -}  & F(A)\CAlg. \xtwocell[-1,-1]{}\omit{<1>} }\]
\bprf First, note that there is a strong symmetric monoidal functor $A\otimes -:\V\to \AMod$, which therefore induces the functor at the top of the diagram, and similarly for the one in the bottom.  The isomorphism \[\nabla: F(A)\otimes F(B) \to F(A\otimes B)\] natural in $B\in \CMon(\V)$ is given by the structure isomorphism of $F$. The only thing one needs to check is that $\nabla$ is a map of $F(A)$-commutative algebras, but this is a straightforward verification.
\eprf
\end{lema}

\subsection{Simplicial cyclic bar construction in general}
The results in this section are similar to the ones in \cite[Section 1]{stonek-graded} which are done for the simplicial \emph{reduced} bar construction.

\bdefn The \emph{simplicial cyclic bar construction} is a functor \[B^\cy_\bullet:\Mon(\V) \to s\V\]
defined as follows. If $A\in \Mon(\V)$ with multiplication $\mu:A\otimes A\to A$ and unit $\eta:\1\to A$, then $B^\cy_n(A)=A^{\otimes n+1}$. The faces $d_i:A^{\otimes n+1}\to A^{\otimes n}$, $i=0,\dots,n$ are defined as
\[d_i= \id^{\otimes i} \otimes \mu  \otimes \id^{\otimes n-i-1} \hspace{.5cm} \text{if } i=0,\dots,n-1, \textup{ and}\]
\[d_n=(\mu \otimes \id^{\otimes (n-1)}) \circ \sigma_{n+1}\]
where $\sigma_{n+1}:A^{\otimes n+1}\to A^{\otimes n+1}$ is the isomorphism that puts the last $A$ term at the beginning. The degeneracies $s_i:A^{\otimes n+1}\to A^{\otimes n+2}$ are \[s_i=\id^{\otimes i+1} \otimes \eta \otimes \id^{\otimes n-i} \hspace{.5cm} \text{ for all }i=0,\dots,n.\]
\edefn
This is a strong symmetric monoidal functor, so it induces a functor between the categories of commutative monoids. But a version of the Eckmann-Hilton argument (see e.g. \cite[6.29]{aguiar}) says that $\CMon(\Mon(\V))\cong \CMon(\V)$, so we get a functor
\[B_\bullet^\cy:\CMon(\V)\to s\CMon(\V).\]

For a commutative monoid $A$ in $\V$, we have that $B_\bullet^\cy(A)\in s\ACAlg$. Indeed, there is a natural morphism $cA\to B_\bullet^\cy A$ in $s\CMon(\V)$, where $cA$ denotes the constant simplicial object at $A$. In simplicial level $n$, it is $\id \otimes \eta^{\otimes n}:A\to A^{\otimes n+1}$. We could specify the simplicial commutative $A$-algebra structure of $B^\cy_\bullet(A)$ more explicitely: the $A$-module structure on $A^{\otimes n+1}\cong A\otimes A^{\otimes n}$ is given by acting on the first factor, and the multiplication over $A$ is given by 
\[\xymatrix{A^{\otimes n+1}\otimes_A A^{\otimes n+1}\cong A \otimes (A^{\otimes n} \otimes A^{\otimes n}) \ar[r]^-{\id \otimes \mu} & A\otimes A^{\otimes n}}\]
where $\mu$ denotes the product of $A^{\otimes n}\in \CMon(\V)$.

Note moreover that $B_\bullet^\cy(A)$ admits a map of simplicial commutative monoids to $cA$: in simplicial level $n$, it is the multiplication of $n$ elements of $A$, which makes sense by commutativity of $A$. So $B_\bullet^\cy(A)$ is a simplicial augmented commutative $A$-algebra.\\

There is a relative version of this construction: if $M$ is an $A$-bimodule, then one can define $B_\bullet^\cy(A,M)$ which has $B_n^\cy(A,M)=M\otimes A^{\otimes n}$ with similar faces and degeneracies. If $M$ is a commutative $A$-algebra, then $B_\bullet^\cy(A,M)$ is a simplicial augmented commutative $M$-algebra.\\

Let $F:\V\to \W$ be a strong symmetric monoidal functor between cocomplete closed symmetric monoidal categories which preserves colimits. Since it takes commutative $A$-algebras to commutative $F(A)$-algebras, then $F(B^\cy_\bullet A)$ is a simplicial commutative $F(A)$-algebra. The structure morphisms of $F$ provide an %
 isomorphism \begin{equation} \label{bf} B_\bullet^\cy(FA)\stackrel{\cong}{\to} F(B_\bullet^\cy A)\end{equation} of simplicial augmented commutative $F(A)$-algebras.

\subsection{Geometric realization} \label{sect:geom} Consider $F=\sip:\Top\to \SMod$: it is a strong symmetric monoidal left adjoint functor \cite[II.1.2]{ekmm}. If $G$ is a topological commutative monoid, we denote by $\S[G]$ the $\S$-module $\sip G$ together with the commutative $\S$-algebra structure induced by the monoid structure of $G$. Note that the map $G\to *$ gives $\S[G]$ an augmentation $\S[G]\to \S$.

Endow the cartesian category $\Top$ with the standard cosimplicial space $\Delta_\top^\bullet$ and the symmetric monoidal category $\SMod$ with the cosimplicial spectrum $\sip \Delta^\bullet_\top$. By \cite[2.9]{stonek-graded}, 
these beget strong symmetric monoidal functors of geometric realization
\[|-|:s\Top \to \Top \hspace{1cm} \textup{and} \hspace{1cm} |-|:s\SMod \to \SMod.\]

If $A$ is a topological commutative monoid or a commutative $\S$-algebra, define \[B^\cy(A)\coloneqq|B^\cy_\bullet(A)|.\]
It is an augmented commutative $A$-algebra. In the $\S$-module case, this object defines the topological Hochschild homology of $A$, denoted $THH(A)$ (which coincides with the derived smash product $A \wedge^L_{A^e} A$ when $A$ is a cofibrant commutative $\S$-algebra \cite[IX.2.7]{ekmm}), %
and if $M$ is an $A$-bimodule or commutative $A$-algebra then $|B_\bullet^\cy(A,M)|$ defines $THH(A,M)$. Note that if $f:A\to B$ is a weak equivalence of cofibrant commutative $\S$-algebras, then $THH(A)\to THH(B)$ is a weak equivalence. First note that $B_\bullet^\cy(A)\to B_\bullet^\cy(B)$ is a weak equivalence in each level. Indeed, $f^{\wedge p+1}$ is a weak equivalence since $f$ is a weak equivalence and $A$ and $B$ are cofibrant. %
Then, since the simplicial cyclic bar construction gives a proper simplicial $\S$-module \cite[IX.2.8]{ekmm}, we can apply \cite[X.2.4]{ekmm} to get the conclusion. %

Since the functors $|-|$ are strong symmetric monoidal, by realizing the isomorphism (\ref{bf}) we obtain the following 

\bprop \label{conmutarTHH} Let $G$ be a topological commutative monoid. There is an isomorphism of augmented commutative $\S[G]$-algebras
\[\xymatrix{THH(\S[G])\ar[r]^-\cong & \S[B^\cy(G)].} \]
\eprop

Versions of the previous proposition have appeared in \cite[Remark 4.4]{svw}, in \cite[Theorem 7.1]{hesselholt-madsen97} in the setting of functors with smash product, and in \cite[Example 4.2.2.7]{dgm} in the setting of $\Gamma$-spaces. Note that we take care to note that this isomorphism respects the commutative $\S[G]$-algebra structures.\\

The following proposition is obtained by applying Lemma \ref{lemext} to the functor $\sip$; the isomorphism commutes with the augmentations by inspection.

\bprop \label{Smono} Consider $G, H\in \CMon(\Top)$. There is an isomorphism of augmented commutative $\S[G]$-algebras
\[\xymatrix{\S[G] \wedge \S[H] \ar[r]^-\cong &\S[G \times H]} \]
natural in $H$.
\eprop

\subsection{Cyclic bar construction of a topological abelian group}  \label{section:cyclictop}
Let $G$ be a topological abelian group with unit $0$. %
Denote by $BG$ the model for the classifying space of $G$ which is given by the geometric realization of the reduced bar construction $B_\bullet(0,G,0)$ of $G$: it is a topological abelian group. Moreover, if $G$ is a CW-complex with $0$ a 0-cell and addition is a cellular map, then the same can be said of $BG$. All of this is due to Milgram \cite{milgram}. %

The space $G\times BG$ gets the structure of a commutative $G$-algebra, via the inclusion of the first factor $G\to G\times BG$, which is a morphism of topological abelian groups. It has an augmentation given by projection to the first factor.

\bprop \label{BcyG} Let $G$ be a topological abelian group. There is a homeomorphism of augmented commutative $G$-algebras
\[B^\cy G \cong G \times BG.\]
\bprf %
Let $G_\bullet$ denote the constant simplicial commutative $G$-algebra on $G$. Consider the maps $r_\bullet:B^\cy_\bullet G \to G_\bullet$, $(g_0,\dots,g_p)\mapsto g_0+\dots+g_p$, and $p_\bullet:B^\cy_\bullet G \to B_\bullet G$, $(g_0,\dots,g_p)\mapsto (g_1,\dots,g_p)$. They assemble to a map
\[\xymatrix{B^\cy_\bullet G \ar[r]^-{(r_\bullet, p_\bullet)} & G_\bullet \times B_\bullet G, \hspace{.5cm} (g_0,\dots,g_p)\mapsto (g_0+\dots+g_p,g_1,\dots,g_p).}\]
We also have maps $i_\bullet:G_\bullet \to B_\bullet^\cy G$, $g\mapsto (g,0,\dots,0)$ and $s_\bullet:B_\bullet G \to B^\cy_\bullet G$, $(g_1,\dots,g_p)\mapsto (-g_1-\dots-g_p,g_1,\dots,g_p)$. We sum them up to a map
\[\xymatrix{ G_\bullet \times B_\bullet G \ar[r]^-{i_\bullet+s_\bullet} & B_\bullet^\cy G, \hspace{.5cm} (g,g_1,\dots,g_p)\mapsto (g-g_1-\dots-g_p,g_1,\dots,g_p).}\]
The maps $(r_\bullet,p_\bullet)$ and $i_\bullet+s_\bullet$ are morphisms of simplicial augmented commutative $G$-algebras which are inverse to one another. (Note that the obvious isomorphisms $G\times G^p \cong G^{p+1}$ are not good, because they do not commute with the last face map.) Applying geometric realization we obtain the result.
\eprf
\eprop

A classical result (which we will not use) states that $B^\cy G$ is homotopy equivalent to the free loop space of $BG$ (see e.g. \cite[Section 2]{bhm}).\index{Free loop space}

\subsection{Inverting an element in \texorpdfstring{$THH$}{THH}} %

Let $R$ be a cofibrant commutative $\S$-algebra and $x\in \pi_*R$. Denote by $\eta:R\to THH(R)$ the unit. Since $THH(R)=B^\cy(R)$ is cofibrant as a commutative $\S$-algebra \cite[Lemma 3.6]{svw}, Proposition \ref{base-loc} gives a weak equivalence of augmented commutative $R[x^{-1}]$-algebras %
\begin{equation}\label{comienzo} \xymatrix{THH(R,R[x^{-1}]) \cong R[x^{-1}]\wedge_R THH(R) \ar[r]^-\sim  & THH(R)[\pi_*\eta(x)^{-1}].}\end{equation}
For simplicity, we denote the codomain of this arrow by $THH(R)[x^{-1}]$.

We now aim to prove that $THH(R,R[x^{-1}])$ and $THH(R[x^{-1}])$ are weakly equivalent commutative $R[x^{-1}]$-algebras.  We will obtain this as a consequence of the following more general theorem, by taking the sequence (\ref{seqc}) to be $\S\to R\to R[x^{-1}]$.
\bteo \label{thh-base} Let \begin{equation}\label{seqc}\S\to A\stackrel{f}{\to} B\end{equation} be a sequence of cofibrations of commutative $\S$-algebras. Suppose that the multiplication map $\mu:B\wedge_A B\to B$ is a weak equivalence. Then the map of augmented commutative $B$-algebras
\begin{equation}\label{thh-base-eq}\xymatrix{ B\wedge_A THH(A)\cong THH(A,B) \ar[rr]^-{THH(f,\id)} && THH(B,B)=THH(B)}\end{equation}
is a weak equivalence.
\eteo
This theorem is valid \emph{mutatis mutandis} when $\S$ is replaced by some cofibrant commutative $\S$-algebra.

We draw inspiration from \cite[Lemma 2.4.10]{eva-thesis}. For $R\in \SCAlg$, denote $R^e\coloneqq R\wedge R$. 
\bprf %
Consider $A$ (resp. $B$) as a commutative $A^e$-algebra (resp. $B^e$-algebra) via the multiplication map $A^e\to A$ (resp. $B^e\to B$). %
Recall that $THH(A,B)\cong B\wedge_{A^e} B(A,A,A)$ and similarly for $THH(B)$ (see \cite[IV.7.2, IX.2.4]{ekmm} for a definition of the two-sided bar construction $B(A,A,A)$ and a proof of the isomorphism).

Let $\tilde B \stackrel{\sim}{\to} B$ be a cofibrant replacement of $B$ in the category of commutative $B^e$-algebras. There is a commutative diagram of $\S$-modules %
\begin{equation}\label{sorw}\xymatrix{THH(A,B) \ar[d]_-{THH(f,\id)} & \tilde B \wedge_{A^e} B(A,A,A) \ar[l]_-\sim \ar[r]^-\sim \ar[d]_-{(\id, f)} & \tilde B \wedge_{A^e} A \ar[d]_-{\overline f} \\ THH(B) & \tilde B \wedge_{B^e} B(B,B,B) \ar[r]^-\sim \ar[l]_-\sim & \tilde B \wedge_{B^e} B.}\end{equation}
Indeed, recall that there is a weak equivalence of commutative $A^e$-algebras $B(A,A,A)\to A$ \cite[IV.7.5]{ekmm} and a cofibration in $\SCAlg$ $A^e\to B(A,A,A)$ given by inclusion of the first and last smash factors. See \cite[Proof of Lemma 2.4.8]{eva-thesis} for a proof of this last fact.

The arrow $(\id,f)$ in the middle is defined via the universal property for the coproduct in commutative $A^e$-algebras, using the canonical map $\tilde B \to \tilde B \wedge_{B^e} B(B,B,B)$ to the first factor, and the map $B(A,A,A)\to B(B,B,B)$ defined by smash powers of $f$ at the simplicial level followed by the canonical map to the second factor.

The arrow $\overline{f}$ is described as follows. First note that there are isomorphisms
\[\tilde B \wedge_{A^e} A \cong \tilde B \wedge_{B^e} (B^e \wedge_{A^e} A)\cong \tilde B \wedge_{B^e} (B\wedge_A B).\]
The last step comes from the isomorphism of commutative $B^e$-algebras $B^e \wedge_{A^e} A  \cong B\wedge_A B$
which appears e.g. in \cite[Lemma 2.1]{lindenstrauss}. %
Then $\overline{f}$ is defined to be the composition
\[\xymatrix{\tilde B \wedge_{A^e} A \cong \tilde B \wedge_{B^e} (B\wedge_A B) \ar[r]^-{\id\wedge \mu} & \tilde B \wedge_{B^e} B.}\]
The previous diagram appears as the geometric realization of a diagram in simplicial $\S$-modules. The arrows in this latter diagram are very explicitely defined, and it is immediate that they make the diagram commute.

Therefore, to see that $THH(f,\id)$ is a weak equivalence, it suffices to see that $\id\wedge \mu$ is a weak equivalence. This is the case: indeed, the functor $\tilde B \wedge_{B^e}-$ preserves weak equivalences between cofibrant commutative $\S$-algebras because $\tilde B$ is cofibrant as a commutative $B^e$-algebra. Now note that $B\wedge_A B$ is a cofibrant commutative $\S$-algebra because it is a cofibrant commutative $A$-algebra (it is a coproduct of two cofibrant commutative $A$-algebras).
\eprf

Lemma \ref{idem} allows us to apply Theorem \ref{thh-base} to $\S\to R\to R[x^{-1}]$. Putting this together with the weak equivalence (\ref{comienzo}), we obtain:

\bcor \label{cor-thhloc} Let $R$ be a cofibrant commutative $\S$-algebra, and let $x\in \pi_*R$. There are weak equivalences of augmented commutative $R[x^{-1}]$-algebras
\[\xymatrix{THH(R)[x^{-1}] & THH(R,R[x^{-1}]) \ar[l]_-\sim \ar[r]^-\sim & THH(R[x^{-1}]). }\]
\ecor

\bobs We have recently become aware that, in \cite[Page 353]{svw}, the authors state that ``one can prove that $THH$ commutes with localizations'', pointing to an article in preparation which never appeared.
\eobs

\bobs We know three proofs of the fact that Hochschild homology commutes with localizations. Weibel \cite[9.1.8(3)]{weibel} proves it using the fact that Tor behaves well under flat base change. Brylinski \cite{brylinski} (see also \cite[1.1.17]{loday}) proves it by comparing the homological functors defined on $A$-bimodules $S^{-1}HH_n(A,-)$ and $HH_n(S^{-1}A,S^{-1}-)$, where $S$ is a multiplicative subset of the commutative algebra $A$. In \cite{weibel-geller}, Geller and Weibel prove the more general result that Hochschild homology behaves well with respect to étale maps of commutative algebras $A\to B$, of which a localization map is an example. Our proof of Theorem \ref{thh-base} is closer to the first of these approaches. %
\eobs

\bobs For a map $f:A\to B$ of commutative $\S$-algebras as in Theorem \ref{thh-base}, the question of under what conditions is (\ref{thh-base-eq}) a weak equivalence has been considered before. For example, in \cite[Lemma 5.7]{mccarthy-minasian} the authors prove that it holds when $A$ and $B$ are connective and the unit $B\to THH^A(B)$ is a weak equivalence. Mathew \cite[Theorem 1.3]{mathew-thh}, working in the context of the $E_\infty$-rings of Lurie, proved that a map $A\to B$ of $E_\infty$-rings satisfies that (\ref{thh-base-eq}) is an equivalence provided $f$ is étale, with no hypotheses on connectivity. There is a notion of localization of $E_\infty$-rings, and Lurie proved that localization maps are étale \cite[7.5.1.13]{ha}. This gives a short proof of Theorem \ref{thh-base} applied to $\S\to R\to R[x^{-1}]$ in the context of $E_\infty$-rings: this is the point of view adopted for the more general result of \cite[Corollary 7.4]{rsv-thom}.
\eobs

\section{Topological Hochschild homology of \texorpdfstring{$KU$}{KU}} \label{sect:thhku}

\subsection{Flatness} \label{sect:cofibrancy}

Let $G$ be a topological commutative monoid. We cannot prove that the commutative $\S$-algebra $\S[G]$ is cofibrant (and we believe it is not, in general), even if $G$ satisfies good cofibrancy hypotheses. Instead, we remark that, when $G$ is a CW-complex and the unit is a 0-cell, $\S[G]$ is flat (to be defined below) and $B^\cy_\bullet(\S[G])$ is a proper simplicial $\S$-module: these properties ensure that $THH(\S[G])$ has homotopical meaning. We thank Cary Malkiewich and Michael Mandell for helping us realize this.

\bdefn \label{def-flat}An $\S$-module $M$ is \emph{flat} if $M\wedge-:\SMod\to \SMod$ preserves all weak equivalences.
\edefn

Note that if $M$ is flat and $N$ is any $\S$-module, then $N\wedge M$ computes the derived smash product.

\bobs Flat $\S$-modules satisfy the following properties: \label{pseudo-rem} 
\bi 
\item Any cofibrant $\S$-module is flat, and $\S$ is flat.
\item Smash products of flat $\S$-modules are flat.
\item Coproducts of flat $\S$-modules are flat.
\item Weak equivalences between flat $\S$-modules are closed under finite smash products. %
\ei
We do not know whether the underlying $\S$-module to a cofibrant commutative $\S$-algebra is automatically flat.
\eobs

\begin{lema} \label{cofcofsemi} Let $A$ and $A'$ be cofibrant commutative $\S$-algebras, $N$ be any $\S$-module and $M$ be a flat $\S$-module. Let $f:A\to N$, $g:A'\to M$ be weak equivalences of $\S$-modules. Then 
\[f \wedge g : A \wedge A' \to N \wedge M\]
is a weak equivalence of $\S$-modules.
\bprf Let $\gamma_A:\Gamma A \to A$ and $\gamma_{A'}:\Gamma A' \to A'$ be cofibrant replacements of $A$ and $A'$ in $\SMod$. Consider the following commutative diagram:
\[\xymatrix@C+1pc{A \wedge A' \ar[r]^-{f \wedge g} & N \wedge M \\ \Gamma A \wedge \Gamma A' \ar[u]^-{\gamma_A \wedge \gamma_{A'}} \ar[r]_-{\id \wedge (g\circ \gamma_{A'})} & \Gamma A \wedge M. \ar[u]_-{(f\circ \gamma_A) \wedge \id} }\]
By Properties (\ref{ujin}) and (\ref{es}) of Section \ref{sect-cof}, the left vertical map and the bottom horizontal map are weak equivalences. The right vertical map is a weak equivalence because $M$ is flat. In conclusion, $f\wedge g$ is a weak equivalence.
\eprf
\end{lema}

\begin{lema} \label{isflat} Let $Y$ be a based CW-complex. Then $\Sigma^\infty Y$ is a flat $\S$-module. Also, the functor $\Sigma^\infty Y \wedge -:\S\Mod\to \S\Mod $ is left Quillen,  %
so if $M$ is a cofibrant $\S$-module, then $\Sigma^\infty Y \wedge M$ is so too. %
\bprf The first statement is \cite[4.11(i)]{mandell-may}. Recall that $\Sigma^\infty Y \wedge-$ is isomorphic to $Y\wedge -$, where the latter $\wedge$ denotes the tensor of $\SMod$ over $\Top_*$ \cite[II.1.4]{ekmm}. Since $\SMod$ is a $\Top_*$-model category and $Y$ is cofibrant in $\Top_*$, this implies that $\Sigma^\infty Y \wedge -$ is left Quillen, so it takes the cofibrant $M$ to a cofibrant $\S$-module. %
\eprf
\end{lema}

\bprop \label{pseudocof-thh} Let $G$ be a CW-complex which is a topological commutative monoid with unit a 0-cell. Let $f:\widetilde{\S[G]}\to \S[G]$ be a weak equivalence in $\SCAlg$, where $\widetilde{\S[G]}$ is a cofibrant commutative $\S$-algebra. Then \[THH(f):THH(\widetilde{\S[G]})\to THH(\S[G])\] is a weak equivalence of commutative $\widetilde{\S[G]}$-algebras. %
\bprf First, note that $B_\bullet^\cy(\S[G])$ is a proper simplicial $\S$-module. Indeed, by Proposition \ref{conmutarTHH}, $B_\bullet^\cy(\S[G])\cong \S[B_\bullet^\cy G]$. Now, the functor $\sip:\Top\to \SMod$ preserves properness of simplicial objects, as was observed in \cite[IV.7.8]{ekmm}. But $B^\cy_\bullet(G)$ is a proper simplicial space, since it is good \cite[3.2]{thh-thom} %
and %
any good simplicial space is proper \cite[Proof of 2.4(b)]{lewis-lillig}. Therefore, by \cite[X.2.4]{ekmm}, it suffices to see that $B^\cy_\bullet(f):B^\cy_\bullet(\widetilde{\S[G]})\to B^\cy_\bullet(\S[G])$ is levelwise a weak equivalence.

The map $f\wedge f: \widetilde{\S[G]} \wedge \widetilde{\S[G]} \to \S[G]\wedge \S[G]$ is a weak equivalence by the previous two lemmas.  For higher smash powers of $f$, the statement is proven by an analogous argument and induction.
\eprf
\eprop

\bobs More generally, we have just proven that if $A\to B$ in $\SCAlg$ is a weak equivalence where $A$ is a cofibrant commutative $\S$-algebra and $B$ is a flat $\S$-module such that $B_\bullet^\cy(B)$ is a proper simplicial $\S$-module, then $THH(A)\to THH(B)$ is a weak equivalence of commutative $A$-algebras. %
\eobs

\subsection{Topological Hochschild homology of \texorpdfstring{$\S[G][x^{-1}]$}{S[G][1/x]}}

Let $G$ be a CW-complex which is a topological abelian group with unit a 0-cell and cellular addition map. As remarked in Section \ref{section:cyclictop}, these assumptions guarantee that $BG$ is again a CW-complex, so that $\S[BG]$ is a flat $\S$-module by Lemma \ref{isflat}. This will be useful in the proof of Theorem  \ref{teoloc}. %

We first isolate a result that does not involve inverting an element.%
\begin{lema} \label{sinloc} There is an isomorphism of augmented commutative $\S[G]$-algebras
\[THH(\S[G]) \cong \S[G] \wedge \S[BG] = \S[G][BG].\]
\bprf It is an application of Propositions \ref{conmutarTHH}, \ref{BcyG} and \ref{Smono}:
\[\xymatrix{THH(\S[G]) \ar[r]^-\cong & \S[B^\cy G] \ar[r]^-\cong & \S[G\times BG] & \S[G] \wedge \S[BG] \ar[l]_-\cong.}\qedhere\] 
\eprf
\end{lema}

Let $x\in \pi_* \S[G]$. Recall from Definition \ref{def-nocof} that $\S[G][x^{-1}]_h$ is defined as $(Q\S[G])[x^{-1}]$ where $Q$ is a cofibrant replacement functor in the category of commutative $\S$-algebras coming from a functorial factorization. It factors the unit $e_R:\S\to R$ of a commutative $\S$-algebra $R$ as
\begin{equation}\label{unitrep1} \xymatrix{\S \ar[rr]^-{e_R} \ar@{>->}[rd]_-{e_{QR}} && R \\ & QR. \ar@{->>}[ru]^-\sim_-{q_R}}\end{equation}

\bteo \label{teoloc} The commutative $\S[G][x^{-1}]_h$-algebras $THH(\S[G][x^{-1}]_h)$ and $\S[G][x^{-1}]_h[BG]$ are weakly equivalent as $\S[G][x^{-1}]_h$-algebras.
\eteo

For any commutative $\S$-algebra $A$, the notation $A[BG]$ stands for the commutative $A$-algebra $A\wedge \S[BG]$: thus, its underlying $A$-module is $A\wedge (BG)_+$. %
No confusion should arise from the usage of square brackets for two different notions.

\bprf For ease of notation, let us denote $\S[G]$ by $A$ and $\S[BG]$ by $B$. As in Remark \ref{funct-fact} for $R=QA$, the functor $Q$ begets a cofibrant replacement functor $Q_{QA}$ in the category of commutative $QA$-algebras. For ease of notation we denote it by $Q_A$. This functor factors the unit $u_X:QA\to X$ of a $QA$-commutative algebra $X$ as
\begin{equation} \label{unitrep2}\xymatrix{QA \ar[rr]^-{u_X} \ar@{>->}[rd]_-{u_{Q_AX}} && X \\ & Q_A X. \ar@{->>}[ru]^-\sim_-{q^A_X}}\end{equation}

Using Corollary \ref{cor-thhloc}, we obtain a zig-zag of two weak equivalences of commutative $QA[x^{-1}]$-algebras
\[THH(QA[x^{-1}])\simeq THH(QA)[x^{-1}].\]

By Proposition \ref{pseudocof-thh}, the map of commutative $QA$-algebras $THH(QA)\to THH(A)$ is a weak equivalence. We apply $Q_A$ and obtain a weak equivalence of commutative $QA$-algebras \[\xymatrix{Q_A(THH(QA)) \ar[r]^-\sim & Q_A THH(A).}\]
Note that if we had applied $Q$ instead of $Q_A$, we would not be able to guarantee that the above morphism be a morphism of $QA$-algebras.

We also have a weak equivalence of commutative $QA$-algebras $Q_A THH(QA) \to THH(QA)$. After inverting $x$, we obtain a zig-zag of weak equivalences of commutative $QA[x^{-1}]$-algebras:
\begin{equation}\label{minu}\xymatrix{THH(QA)[x^{-1}] & (Q_ATHH(QA))[x^{-1}] \ar[l]_-\sim \ar[r]^-\sim & (Q_A THH(A))[x^{-1}].}\end{equation}

Now, by Lemma \ref{sinloc}, $THH(A)\cong A \wedge B$ as commutative $A$-algebras, so we obtain an isomorphism of commutative $QA[x^{-1}]$-algebras
\[(Q_A THH(A))[x^{-1}] \cong (Q_A(A \wedge B))[(x\wedge 1)^{-1}].\]

We now construct a weak equivalence of commutative $QA$-algebras \[ r: Q_AA \wedge QB \to Q_A(A \wedge B)\]
where $Q_AA\wedge QB$ is a $QA$-algebra by means of the map $\xymatrix{QA\ar[r]^-{u_{Q_AA}} & Q_AA \ar[r]^-{\id \wedge e_{QB}} & Q_AA \wedge QB.}$ Since $\wedge$ is the coproduct in the category of commutative $\S$-algebras, we define the map $r$ to be the morphism in $\SCAlg$ given as follows. It is $Q_A(\id \wedge e_B):Q_AA \to Q_A(A\wedge B)$ on the first component. On the second component, it is the map $t:QB\to Q_A(A\wedge B)$ gotten from the functorial factorization on $\SCAlg$ applied to the vertical arrows in the following commutative diagram:
\[\xymatrix@C+1pc{\S \ar[r]^-{e_{QA}} \ar[d]_-{e_B} & QA \ar[d]^-{q_A \wedge e_B} \\ B\ar[r]_-{e_A\wedge \id_B} & A \wedge B.  }\]
This defines the map $r$ as a morphism of commutative $\S$-algebras. It is a morphism of commutative $QA$-algebras: the composition $\xymatrix@C+1pc{QA \ar[r]^-{u_{Q_AA}} & Q_AA \ar[r]^-{Q_A(\id\wedge e_B)} & Q_A(A\wedge B)}$ coincides with $u_{Q_A(A\wedge B)}$ by functoriality of the factorization. We will now prove that $r$ is a weak equivalence. First, consider the following diagram:
\[\xymatrix@C+2pc{Q_AA \ar[r] \ar[rd]_-{Q_A(\id\wedge e_B)} \ar@/_2pc/[rdd]_-{q^A_A\wedge e_B} & Q_AA \wedge QB \ar[d]^-r & QB \ar[l] \ar[ld]^-t \ar@/^2pc/[ldd]_-{\ \ e_A\wedge q_B} \\ & Q_A(A\wedge B) \ar[d]^-{q^A_{A\wedge B}} \\ & A\wedge B  }\]
The two upper triangles commute by definition of $r$. The two lower triangles also commute: the left one by definition of $Q_A(\id\wedge e_B)$, the right one by definition of $t$. By the universal property of the coproduct in $\SCAlg$, this proves that $q^A_{A\wedge B} \circ r=q_A^A \wedge q_B$. %
Lemma \ref{cofcofsemi} applies to prove that $q_A^A \wedge q_B$ is a weak equivalence. Therefore, since $q^A_{A\wedge B}$ is also a weak equivalence, so is $r$.

Inverting $x\wedge 1$ in $r$, we obtain a weak equivalence of commutative $QA[x^{-1}]$-algebras
\[ %
\xymatrix{Q_A(A \wedge B)[(x\wedge 1)^{-1}] & \ar[l]_-\sim (Q_AA \wedge Q B)[(x\wedge 1)^{-1}].}\]
By Proposition \ref{sepa}, there is a weak equivalence of commutative $QA[x^{-1}]$-algebras
\begin{equation}\label{vers}\xymatrix{(Q_AA \wedge Q B)[(x\wedge 1)^{-1}] & \ar[l]_-\sim Q_AA[x^{-1}] \wedge QB. }\end{equation}
From (\ref{unitrep2}) applied to $X=A$, we get that $q_A=q_A^A\circ u_{Q_AA}$, so $u_{Q_AA}:QA\to Q_AA$ is a weak equivalence of commutative $QA$-algebras. %
We can invert $x$ to obtain the weak equivalence of commutative $QA[x^{-1}]$-algebras $QA[x^{-1}]\to Q_AA[x^{-1}]$, which after smashing with $QB$ becomes
\[\xymatrix{Q_A A[x^{-1}] \wedge QB & QA[x^{-1}] \wedge QB, \ar[l]_-\sim }\] %
a weak equivalence of commutative $QA[x^{-1}]$-algebras.

Now consider the weak equivalence $q_B:QB\to B$ of commutative $\S$-algebras. Lemma \ref{cofcofsemi} applies to prove that 
\[\xymatrix{QA[x^{-1}]\wedge QB \ar[r]^-\sim & QA[x^{-1}] \wedge B }\]
is a weak equivalence of commutative $QA[x^{-1}]$-algebras.
Putting together all these weak equivalences, we obtain a zig-zag of weak equivalences of commutative $(Q\S[G])[x^{-1}]$-algebras \linebreak $THH((Q\S[G])[x^{-1}])\simeq ((Q\S[G])[x^{-1}])[BG]$.
\eprf

\bobs \label{augment-complicado} In this remark, we explain how compatible is the zig-zag of weak equivalences just obtained with the augmentations.

Ignoring the last weak equivalence in the zig-zag of weak equivalences obtained in the proof, we have obtained
\begin{equation}\label{zaq} THH(QA[x^{-1}])\simeq QA[x^{-1}]\wedge QB. \end{equation}
The left-hand side is naturally augmented over $QA[x^{-1}]$. We now give the right-hand side an augmentation over $QA[x^{-1}]$ defined from the universal property of the coproduct in $\SCAlg$. On the first factor, it is the identity. On the second factor, it is the arrow:
\[\xymatrix@C+1pc{QB \ar[r]^-{Q(e_A \circ \epsilon_B)} & QA \ar[r] & QA[x^{-1}] }\]
where $\epsilon_B$ is the augmentation of $B=\S[BG]$ coming from $BG\to *$ and the second map is the localization map.
If one goes through the proof of the previous theorem, one can prove that the zig-zag (\ref{zaq}) commutes with the augmentations, in the following sense: this zig-zag fits as one long horizontal side of a ladder diagram, where the other side starts and ends with $QA[x^{-1}]$ and in the middle has as objects $QA[x^{-1}]$, $Q_AQA[x^{-1}]$ or $Q_AA[x^{-1}]$, connected to each other via the obvious weak equivalences or via identities when possible. The ingredients used in the proof of this are: the fact that Corollary \ref{cor-thhloc} and Lemma \ref{sinloc} are compatible with the augmentations, the naturality of Proposition \ref{sepa} and of $THH$, and the functoriality of the factorization (cofibrations, acyclic fibrations).

The last weak equivalence
\begin{equation}\label{xlo} QA[x^{-1}]\wedge QB \to QA[x^{-1}]\wedge B\end{equation}
is more problematic. The only sensible augmentation over $QA[x^{-1}]$ to define on the codomain would be $\id \wedge \epsilon_B$, but then the weak equivalence (\ref{xlo}) commutes with the augmentations only up to homotopy of commutative $QA[x^{-1}]$-algebras. Indeed, in the following diagram
\[\xymatrix@C+2pc{QB\ar[r]^-{Q(e_A \circ \epsilon_B)} \ar[d]_-{q_B} & QA \ar[d]^-{q_A}_-\sim \\ B \ar[ru]^-{e_{QA} \circ \epsilon_B} \ar[r]_-{e_A \circ \epsilon_B} & A}\]
the square commutes, the bottom triangle commutes, but the upper triangle does not seem to commute. However, since $q_A$ is a weak equivalence, it commutes in the homotopy category of $QA\CAlg$.
\eobs

\subsection{Snaith's theorem and first description of $THH(KU)$} \label{construction-ku} \index{KU} \index{Topological $K$-theory}

There is a cofibrant commutative $\S$-algebra $KU$ of complex topological $K$-theory \cite[VIII.4.3]{ekmm}. It is obtained by applying the localization theorem we reviewed in Theorem \ref{rect} to the cofibrant commutative $\S$-algebra $ku$ of connective complex $K$-theory and its Bott element. Here $ku$ is constructed by multiplicative infinite loop space theory.

The presentation for $KU$ which we will use relies on the following theorem of Snaith \cite{snaith79}, \cite{snaith81}. %

\bteo $KU$ is weakly equivalent as a homotopy commutative ring spectrum to \linebreak $\S[\mathbb C P^\infty][x^{-1}]_{\textup{tel}}$, where $x\in \pi_2(\S[\C P^\infty])$ is represented by the map induced from the inclusion $\C P^1\to \C P^\infty,$ i.e.
\begin{equation}\label{xproj} \Sigma^\infty S^2\cong \Sigma^\infty \C P^1 \to \S \vee \Sigma^\infty \C P^\infty \simeq \sip \C P^\infty. \end{equation}
Here $\S[\C P^\infty][x^{-1}]_{\textup{tel}}$ means the homotopy commutative ring spectrum obtained with a telescope construction, as in (\ref{telescope}).
\eteo

As remarked in \cite[VIII.4]{ekmm}, the more structured version of the inversion of a homotopy element described in Section \ref{inversion} is weakly equivalent to this telescope construction as homotopy commutative ring spectra (i.e. commutative monoids in the homotopy category of spectra). %
Indeed, the technology of $\S$-algebras did not exist at the time Snaith's theorem got published, but this is not a problem:

\bteo \label{bari} \cite[6.2]{baker-richter} Let $A$ be an $E_\infty$-ring spectrum which is weakly equivalent to $KU$ as a homotopy commutative ring spectrum. Then $A$ is weakly equivalent to $KU$ as an $E_\infty$-ring spectrum.
\eteo

Any $E_\infty$-ring spectrum $A$ admits a weak equivalence of $E_\infty$-ring spectra from the commutative $\S$-algebra $\S\wedge_\L A$, and this construction is functorial \cite[II.3.6]{ekmm}; moreover, the morphisms of commutative $\S$-algebras are exactly the morphisms of the underlying $E_\infty$-ring spectra. So if $A$ in the statement of Theorem \ref{bari} was a commutative $\S$-algebra to begin with, then it is weakly equivalent to $KU$ as a commutative $\S$-algebra. %

Let $K(\Z,2)$ denote the topological abelian group given by $B(B\Z)$, where $B$ is the classifying space construction reviewed in Section \ref{section:cyclictop}.  
The homotopy commutative ring spectrum $\S[\C P^\infty]$ is weakly equivalent to the cofibrant commutative $\S$-algebra $Q\S[K(\Z,2)]$. %
Therefore, \[\S[\C P^\infty][x^{-1}]_{\textup{tel}} \simeq Q\S[K(\Z,2)][x^{-1}]\]
as homotopy commutative ring spectra. By the results above, we obtain that
\begin{equation}\label{defku} KU \simeq Q\S[K(\Z,2)][x^{-1}]\end{equation}
as commutative $\S$-algebras. This is the description of $KU$ that we shall be using, so from now on we let $KU$ denote the cofibrant commutative $\S$-algebra $Q\S[K(\Z,2)][x^{-1}]$ that we have also denoted by $\S[K(\Z,2)][x^{-1}]_h$.

\bobs We thank Christian Schlichtkrull for pointing out the article \cite{arthan} to us. In Theorems 5.1 and 5.2 therein, it is proven that if $t\in \pi_n(\S[K(\Z,n)])$ is a generator, %
then the spectrum $\S[K(\Z,n)][t^{-1}]_{\textup{tel}}$ is contractible for $n$ odd and is equivalent to $H\Q[t^{\pm 1}]$ for $n\geq 4$ even. So the case $n=2$ which we treat here is the only interesting localization of $\S[K(\Z,n)]$.
\eobs

Note that $\Z$ is a CW-complex (with only 0-cells) and a (discrete) topological abelian group with cellular addition map, so this guarantees that $K(\Z,2)$ satisfies the same hypotheses, as recalled in Section \ref{section:cyclictop}. Therefore, we can apply Theorem \ref{teoloc} to obtain:

\bteo \label{thhku1} The commutative $KU$-algebras $THH(KU)$ and $KU[K(\Z,3)]$ are weakly equivalent as commutative $KU$-algebras.
\eteo

Remark \ref{augment-complicado} tells us that the zig-zag of weak equivalences is compatible up to homotopy with the augmentations.

\bobs \label{thh-thom} Compare with what happens to $THH(MU)$: in \cite{thh-thom}, the authors establish a weak equivalence of $\S$-modules $THH(MU)\simeq MU \wedge SU_+$. They actually prove the following more general result. Let $BF$ denote a classifying space for stable spherical fibrations. %
If $f:X\to BF$ is a 3-fold loop map and $T(f)$ is its Thom spectrum, then there is a weak equivalence of $\S$-modules \begin{equation}\label{thhtf} THH(T(f))\simeq T(f)\wedge BX_+.\end{equation}Note that this result was improved to a weak equivalence of $E_\infty$ $\S$-algebras by Schlichtkrull \cite[Corollary 1.2]{schlichtkrull-higher} in the case where $X$ is a grouplike $E_\infty$-space and $f$ is an $E_\infty$-map.

Our Theorem \ref{thhku1} gives in particular a weak equivalence of $\S$-modules \[THH(KU)\simeq KU \wedge K(\Z,3)_+.\]By comparing this formula to (\ref{thhtf}), one is naturally led to conjecture that $KU$ is the Thom spectrum of an $\infty$-loop map $K(\Z,2)\simeq BU(1)\to BU$. However, this is not possible, since such Thom spectra are connective. In the last decade, more general Thom spectra which can be non-connective have been introduced, and in \cite[Example 4.23]{rsv-thom} the authors remarked that $KU$ cannot be such a Thom spectrum of a map from $K(\Z,2)$, as a consequence of the Thom isomorphism theorem. One possible explanation of why does $KU$ behave like a Thom spectrum to the eyes of topological Hochschild homology is to be found in that paper, where the above expression for $THH(KU)$ is obtained by considering $KU$ as an étale extension of $\S[K(\Z,2)]$, which is a trivial Thom spectrum.
\eobs

\subsection{Rationalization} \index{Rationalization}In this section we review some facts about rationalization of $\S$-modules and of based spaces that we will be using, often without explicit mention.\\

\label{rationalization} %

Consider a model for the Eilenberg-Mac Lane spectrum of $\Q$ which is a commutative $\S$-algebra \cite[II.4]{ekmm}. Let $H\Q$ be a cofibrant $\S$-module equivalent
to it in the category of $\S$-modules. In particular, $H\Q$ is a homotopy commutative ring spectrum. Denote by $\iota:\S\to H\Q$ its unit and by $\mu:H\Q\wedge H\Q\to H\Q$ its multiplication map, which is a weak equivalence.

We consider Bousfield localization \cite[VIII.1]{ekmm}, \cite[19.2]{mayponto} of $\S$-modules with respect to $H\Q$, and we call this process \emph{rationalization}. A map $X\to Y$ of $\S$-modules is an \emph{$H\Q$-equivalence} (or \emph{rational equivalence}) if it is a weak equivalence after smashing it with $H\Q$. An $\S$-module $W$ is \emph{$H\Q$-acyclic} if $H\Q \wedge W\simeq *$. An $\S$-module $X$ is \emph{$H\Q$-local} (or \emph{rational}) if, for every $H\Q$-acyclic $\S$-module $W$, the only map $W\to X$ in the homotopy category of $\S$-modules is the trivial one. A \emph{rationalization map} for $X$ is an $H\Q$-equivalence $X\to Y$ where $Y$ is rational; rationalizations are unique up to homotopy. %
Note that $H\Q\wedge X$ is rational since it is an $H\Q$-module \cite[1.17]{ravenel-loc}. %
The following diagram is homotopy commutative:
\[\xymatrix@C+1pc{H\Q\wedge X \ar[r]^-{\id\wedge \iota \wedge \id} \ar[rd]_-{\id} & H\Q \wedge H\Q \wedge X\ar[d]^-{\mu \wedge \id} \\ & H\Q \wedge X,}\]
so we only need $\mu\wedge \id$ to be a weak equivalence in order to assert that $\iota \wedge \id: X\to H\Q \wedge X$ is a rationalization of $X$. This is true because, as we have seen in Property (\ref{moguil}) of Section \ref{sect-cof}, all $\S$-modules $X$ satisfy that $X\wedge -$ preserves weak equivalences between cofibrant $\S$-modules. Therefore, $\iota \wedge \id: X\to H\Q\wedge X$ is a rationalization map for $X$, and from now on we let $X\to X_\Q$ mean this map. This construction is functorial in $X\in\SMod$. %

From this construction one deduces some properties: if $f:X\to Y$ is a weak equivalence of $\S$-modules then it is a rational equivalence, since $H\Q$ is cofibrant. The homotopy groups of $X_\Q$ are isomorphic to $\Q \otimes \pi_*X$. %
An $\S$-module is rational if and only if $X\to X_\Q$ is a weak equivalence, 
if and only if the homotopy groups of $X$ are rational (i.e. $\Q$-vector spaces). %
If $f:X\to Y$ is a rational equivalence between rational $\S$-modules, then it is a weak equivalence.  %
If $X$ is the suspension spectrum of a based CW-complex, then $X_\Q=H\Q\wedge X$ is cofibrant (Lemma \ref{isflat}).

Let $Y$ be another $\S$-module. 
We have a map \[X_\Q \wedge Y_\Q \to (X\wedge Y)_\Q\] given by applying the multiplication map of $H\Q$, and it is a weak equivalence. 
In particular, the smash product of the rationalization maps, $X\wedge Y \to X_\Q \wedge Y_\Q$, is a rationalization of $X\wedge Y$. Note that if $X$ and $Y$ are rational and one of them is cofibrant or flat, so that $X\wedge Y$ computes the derived smash product,
 then $\pi_*(X \wedge Y)\cong \pi_*X \otimes_\Q \pi_*Y$ by a Künneth spectral sequence argument computing $H\Q_*(X\wedge Y)$ \cite[IV.4.7]{ekmm}. %

Let $n$ be any integer. The degree $n$ map $n:\S\to \S$ induces a map $n:X\to X$ on any $\S$-module $X$ by smashing with it. If $p:X\to X$ is a weak equivalence for every prime $p$ 
then the homotopy groups of $X$ are rational, since $p$ induces the multiplication by $p$ map on homotopy groups. %
Therefore, in this case, $X$ is rational.\\

Based spaces $X$ also admit rationalizations $q:X\to X_\Q$ (see \cite[9.(b)]{fht-rational} or \cite[6.5]{mayponto} for the simply-connected case which is the one we shall be using). %
We will need the following fact concerning the rationalization of integral Eilenberg-Mac Lane spaces \cite[Page 202]{fht-rational}: for $n\geq 2$, there are homotopy equivalences 
\[K(\Z,n)_\Q \stackrel{\sim}{\leftarrow} \begin{cases} S^n_\Q & \textup{if } n \textup{ is odd}, \\ \Omega S^{n+1}_\Q & \textup{if } n \textup{ is even.} \end{cases}\] 
Actually, the authors prove that, for $n$ even, there is a rational homotopy equivalence $\Omega S^{n+1} \to K(\Z,n)$, so we get a homotopy equivalence $(\Omega S^{n+1})_\Q \stackrel{\sim}{\to} K(\Z,n)_\Q$, which is not exactly what we wrote. But more generally, we have that $(\Omega X)_\Q$ is homotopy equivalent to $\Omega X_\Q$. Indeed, %
by taking homotopy groups, we quickly see that $\Omega X_\Q$ is rational and that $\Omega q:\Omega X\to \Omega X_\Q$ is a rational homotopy equivalence. Since rationalizations are unique up to homotopy, this gives the result.\\

If $X$ is a based simply-connected space with rationalization map $q:X\to X_\Q$, then $\Sigma^\infty q:\Sigma^\infty X \to \Sigma^\infty X_\Q$ is immediately seen to be a rationalization of $\Sigma^\infty X$. %
Therefore, $(\Sigma^\infty X)_\Q$ is homotopy equivalent to $\Sigma^\infty X_\Q$, i.e. rationalization of based spaces and of $\S$-modules are compatible under the $\Sigma^\infty$ functor.

\subsection{$THH(KU)$, continuation} \label{sect-thhku2}

We will now describe the commutative $KU$-algebra $THH(KU)$ as the free commutative $KU$-algebra on the $KU$-module $\Sigma KU_\Q$, %
and we will prove this algebra is weakly equivalent to the split square-zero extension of $KU$ by $\Sigma KU_\Q$. Let us first define this concept. 

Let $R$ be a commutative $\S$-algebra, let $A$ be a commutative $R$-algebra and let $M$ be a non-unital commutative $A$-algebra. Then $A\vee M$ (coproduct of $A$-modules) has a commutative $A$-algebra structure. Indeed, after distributing, a multiplication map
\[(A\vee M) \wedge_A (A\vee M) \to A\vee M\]
looks like
\begin{equation}\label{sqz}(A\wedge_A A) \vee (A\wedge_A M) \vee (M\wedge_A A) \vee (M\wedge_A M) \to A\vee M.\end{equation}
We may define a map like (\ref{sqz}) by defining maps from each of the wedge summands to $A\vee M$. %
Define the maps to $A\vee M$ from $A\wedge_A A$, $A\wedge_A M$ and $M\wedge_A A$ to be the canonical isomorphisms followed by the canonical maps into the respective factor. Finally, consider the map $M\wedge_A M\to A\vee M$ given by the multiplication map of $M$ followed by the canonical map to $A\vee M$. We have thus defined a multiplication map (\ref{sqz}) such that $A\vee M$ is a commutative $A$-algebra with unit given by the canonical map $A\to A\vee M$. We say that $A\vee M$ is a \emph{split extension of $A$ by $M$}. Note that $A\vee M$ is augmented over $A$: the augmentation is the identity on $A$ and the trivial map on $M$. 
If the multiplication of $M$ is trivial, then $A\vee M$ is a \emph{split square-zero extension of $A$ by $M$}; in this case, $M$ is no more than an $A$-module. %

The rest of this section is devoted to the proof of the following

\bteo \label{thhku2} There are weak equivalences of commutative $KU$-algebras
\[\xymatrix{KU \vee \Sigma KU_\Q & F(\Sigma KU_\Q) \ar[l]_-h^-\sim \ar[r]^-{\tilde f}_-\sim & THH(KU)}\]
where $KU \vee \Sigma KU_\Q$ is a split square-zero extension. Here $f:\Sigma KU_\Q\to THH(KU)$ is a morphism of $KU$-modules to be constructed in (\ref{defff}), and the morphism $\tilde f$ is induced from $f$ by the free commutative algebra functor $F:KU\Mod\to KU\CAlg$. The map $h$ is adjoint to the wedge inclusion of $KU$-modules $\Sigma KU_\Q \to KU \vee \Sigma KU_\Q$.
\eteo

\bobs \label{f-monadic} The functor $F$, or more generally, the free commutative algebra functor $F_R:\RMod\to R\CAlg$ where $R$ is a commutative $\S$-algebra, is the left adjoint of the forgetful functor $U_R: R\CAlg\to\RMod$, or alternatively, the free algebra functor for the monad $\mathbb P_R$ on $\RMod$ defined as 
\begin{equation}\label{monadwedge} \mathbb P_R(M)=\bigvee_{n\geq 0} \bigslant{M^{\wedge_R n}}{\Sigma_n} = R \vee M \vee \bigvee_{n\geq 2} \bigslant{M^{\wedge_R n}}{\Sigma_n},\end{equation} where $\Sigma_n$ is the symmetric group on $n$ elements (see e.g. \cite[II.7.1]{ekmm} or \cite[Section 1]{basterra}). Note that $F_RM$ is augmented over $R$: the augmentation is the identity on the $0$-th term and the trivial map on the other terms. %

As explained in Section \ref{sect-cof}, the functor $U_R:\RCAlg\to \RMod$ is a right Quillen functor, so $F_R:\RMod \to \RCAlg$ is a left Quillen functor. In particular, it preserves weak equivalences between cofibrant $R$-modules. %

Note as well that, if $R$ is a cofibrant commutative $R$-algebra and $M\in \RMod$ is %
cofibrant, then the arrow $\bigvee_{n\geq 0} (M^{\wedge_R n})_{h\Sigma_n}\to F_R(M)$ induced from the canonical arrows from the homotopy orbits to the orbits is a %
weak equivalence \cite[III.5.1]{ekmm}. This is a step in the proof of the determination of the model structure on $\RCAlg$. %
\eobs

\bobs \label{antecedentes} A spectrum-level result related to Theorem \ref{thhku2} was obtained by McClure and Staffeldt in \cite[Theorem 8.1]{mc-st}: they showed that $THH(L)\simeq L \vee \Sigma L_\Q$ as spectra, where $L$ is the $p$-adic completion of the Adams summand of $KU$ for a given odd prime $p$; the result was extended to $p=2$ by Angeltveit, Hill and Lawson in \cite[2.3]{thhko}. Ausoni \cite[Proposition 7.13]{ausoni-thhku} formulated without proof the analogous theorem (for an odd $p$) for $KU$ completed at $p$ in place of $L$. In Corollary 7.9 of \cite{thhko}, the authors show that $THH(KO)\simeq KO\vee \Sigma KO_\Q$ as $KO$-modules. The methods used in the proofs of the results just cited are different from ours.
\eobs

We first prove some results needed for the proof. Note that in the following statement we are considering $K(\Z,3)$ as a based space: we are not adding a disjoint basepoint.

\begin{prop} \label{kurat} There is a zig-zag of weak equivalences of $KU$-modules \[KU\wedge K(\Z,3) \simeq \Sigma KU_\Q.\]
\bprf Let $p$ be a prime and consider the cofiber sequence of $KU$-modules
\begin{equation}\label{cofib}\xymatrix{KU \ar[r]^-p & KU \ar[r] & KU/p %
.}\end{equation}
If $p>2$, then $KU/p$ is equivalent to $\bigvee\limits_{i=0}^{p-2} \Sigma^{2i} K(1)$ (see \cite[Lecture 4]{adams-lect}), where $K(1)\simeq L/p$ is the first Morava $K$-theory at $p$. If $p=2$, then $K(1) \simeq KU/2$.

The homology $K(1)_* K(\Z,3)$ is trivial: see \cite[Theorem 12.1]{rw80} for the $p>2$ case, and \cite[Appendix]{jw85} for the $p=2$ case. %
Therefore, after smashing (\ref{cofib}) with $K(\Z,3)$, we get a weak equivalence of $KU$-modules \[\xymatrix{KU \wedge K(\Z,3) \ar[r]_-{p \wedge \id}^\sim & KU \wedge K(\Z,3)}\]
for all primes $p$. This means that $KU\wedge K(\Z,3)$ is rational, and so we have weak equivalences
\begin{align*}
\xymatrix{
KU \wedge K(\Z,3) \ar[r]^-\sim &
(KU \wedge K(\Z,3))_\Q \\
&
\ar[u]_-\sim KU_\Q \wedge K(\Z,3)_\Q &
\ar[l]_-\sim KU_\Q \wedge S^3_\Q \ar[r]^-\sim &
(KU \wedge S^3)_\Q \ar[r]^-\sim &
\Sigma KU_\Q}\end{align*}
by the results quoted in Section \ref{rationalization}, plus Bott periodicity for the last step.
\eprf
\end{prop}

\begin{lema} \label{lema-htpyorbits} The $\S$-modules $\bigslant{(\Sigma H\Q)^{\sma n}}{\Sigma_n}$ are weakly contractible for all $n\geq 2$. %
\begin{proof}
First, note that $\Sigma H\Q$ is a cofibrant $\S$-module. Indeed, $H\Q$ is a cofibrant $\S$-module and $S^1$ is a CW-complex, so by Lemma \ref{isflat} $\Sigma H\Q=\Sigma^\infty S^1 \sma H\Q$ is a cofibrant $\S$-module. Therefore, we can apply \cite[III.5.1]{ekmm} to deduce that the map from the homotopy orbits \[((\Sigma H\Q)^{\sma n})_{h\Sigma_n} \to \bigslant{(\Sigma H\Q)^{\sma n}}{\Sigma_n}\] is a weak equivalence. We will prove that the homotopy orbits form a weakly contractible $\S$-module, thus finishing the proof. %

The homotopy orbits spectral sequence \cite[3.2]{tsalidis} here looks like this:
\[E^2_{*,*} \cong H_*(\Sigma_n; \pi_*((\Sigma H\Q)^{\sma n})) \Rightarrow \pi_*(((\Sigma H\Q)^{\sma n})_{h\Sigma_n}).\]
Remark that $\pi_*((\Sigma H\Q)^{\sma n}) \cong \widetilde{H}_*((S^1)^{\sma n}; \Q)$ as $\Sigma_n$-modules, since indeed rearranging the factors in $(\Sigma H\Q)^{\sma n}=(\Sigma^\infty S^1 \sma H\Q)^{\sma n}$ is $\Sigma_n$-equivariant and so is the iterated multiplication map $(H\Q)^{\sma n}\to H\Q$ (at least up to homotopy). %

Since the homology of $\Sigma_n$ with rational coefficients vanishes in positive degrees, %
the only group in the $E^2$-page of the spectral sequence which could be non-trivial is $H_0(\Sigma_n; \widetilde{H}_n((S^1)^{\sma n}; \Q))$. The action of $\Sigma_n$ on $(S^1)^{\sma n}$ permutes the factors. Under the homeomorphism $(S^1)^{\sma n}\cong S^n$, each $\sigma\in \Sigma_n$ acts on $S^n$ by a map $S^n\to S^n$ whose degree is the sign of $\sigma$. Indeed, one can decompose $\sigma$ as a composition of transpositions, which reduces the problem to the determination of the degree of the map $S^2\to S^2$ determined by permuting the two smash factors. This map does have degree $-1$, since it is a reflection of $S^2$ along a plane cutting the sphere in two identical parts. %

In conclusion, the action of $\Sigma_n$ on $\widetilde{H}_n(S^n;\Q)\cong \Q$ is the rational sign representation of $\Sigma_n$, so $H_0(\Sigma_n;\widetilde{H}_n(S^n;\Q))$ vanishes. %
Therefore, the $E^2$-page of the spectral sequence is trivial, and thus $\pi_*(((\Sigma H\Q)^{\sma n})_{h\Sigma_n})=0$, finishing the proof.
\end{proof}
\end{lema}

For the next results, recall that $F_R$ denotes the free commutative $R$-algebra on an $R$-module functor described in Remark \ref{f-monadic}.

\bcor \label{gunb} The map of commutative augmented $\S$-algebras
\[h_\S: F_\S(\Sigma H\Q) \to \S \vee \Sigma H\Q\]
defined as the adjoint to the wedge inclusion of $\S$-modules $\Sigma H\Q \to \S \vee \Sigma H\Q$ is a weak equivalence, 
where $\S\vee \Sigma H\Q$ is a split square-zero extension.
\bprf We have that $F_\S(\Sigma H\Q)= \S \vee \Sigma H\Q \vee \bigvee\limits_{n\geq2} \bigslant{(\Sigma H\Q)^{\wedge n}}{\Sigma_n}$. By construction, $h_\S$ is the identity on the first two wedge summands and it is a trivial map on the $n\geq 2$ summands, so it is a weak equivalence by the previous lemma.
\eprf
\ecor

\begin{prop} \label{split} Let $R$ be a cofibrant commutative $\S$-algebra. %
The map of commutative augmented $\S$-algebras
\[h: F_R(\Sigma R_\Q) \to R \vee \Sigma R_\Q.\]
defined as the adjoint to the wedge inclusion of $R$-modules $\Sigma R_\Q \to R \vee \Sigma R_\Q$ is a weak equivalence, where $R\vee \Sigma R_\Q$ is a split square-zero extension.
\end{prop}

\bobs \label{f-cof}Note that we are applying $F_R$ to a cofibrant $R$-module, and so in particular $F_R(\Sigma R_\Q)$ is a cofibrant commutative $R$-algebra. Indeed, as observed in Lemma \ref{lema-htpyorbits}, $\Sigma H\Q$ is a cofibrant $\S$-module. Now, the extension of scalars functor $R\wedge-: \SMod\to \RMod$ is left Quillen: indeed, its right adjoint, the restriction of scalars functor, is right Quillen since the model structure in $\RMod$ is created through it. Therefore, $R\wedge(\Sigma H\Q)\cong\Sigma R_\Q$ is a cofibrant $R$-module.
\eobs

\bprf Note that for an $\S$-module $X$, we have a natural isomorphism $F_R(R\wedge X)\cong R\wedge F_\S(X)$. Indeed,
\[F_R(R\wedge X)= \bigvee_{n\geq 0} \bigslant{(R\wedge X)^{\wedge_R n}}{\Sigma_n}\cong R\wedge \bigvee_{n\geq 0} \bigslant{X^{\wedge n}}{\Sigma_n} = R\wedge F_\S(X)\]
since the functor $R\wedge-:\S\Mod\to R\Mod$ is a left adjoint and strong symmetric monoidal. %
Therefore, 
\begin{equation}\label{furi} F_R(\Sigma R_\Q)\cong F_R(R\wedge \Sigma H\Q)\cong R\wedge F_\S(\Sigma H\Q).\end{equation}

From Corollary \ref{gunb} we get a weak equivalence of commutative augmented $\S$-algebras $h_\S: F_\S(\Sigma H\Q) \to \S\vee \Sigma H\Q$, and it is readily verified that the following diagram commutes:
\[\xymatrix@C+1pc{R\wedge F_\S(\Sigma H\Q) \ar[r]^-{\id \wedge h_\S} & R \wedge (\S \vee \Sigma H\Q) \ar[d]_-\cong\\ F_R(\Sigma R_\Q) \ar[u]^-\cong \ar[r]_-h & R \vee \Sigma R_\Q.}\] %
Therefore, $h$ is a weak equivalence if and only if $\id \wedge h_\S$ is a weak equivalence, which follows from an application of Lemma \ref{cofcofsemi}.
\eprf

\bobs \label{remark-i} Let $i:R\vee \Sigma R_\Q\to F_R(\Sigma R_\Q)$ be the $R$-module map given by the wedge inclusion. By construction of $h$, we have that $h\circ i=\id$. In particular, $i$ is a weak equivalence. %
\eobs

\bprf[Proof of Theorem \ref{thhku2}] The map $h$ is the one gotten in Proposition \ref{split} for the case $R=KU$. We now aim to establish the equivalence of commutative $KU$-algebras $\tilde f:F(\Sigma KU_\Q)\to THH(KU)$. First, we work additively, and then we will determine the multiplicative structure. 

Recall that for any well-based space $X$, there is a homotopy equivalence of based spaces $\Sigma(X_+)\simeq S^1 \vee \Sigma X$. 
It makes the following diagram of based spaces commute:
\[ \xymatrix{& S^1 \ar[ld] \ar[rd]^-{i_1} \\ \Sigma(X_+) \ar[rr]^-{\simeq} && S^1\vee \Sigma X}\]
where the left diagonal map is $\Sigma u:\Sigma (*_+) \to \Sigma(X_+)$; here $u:*\to X$ is the basepoint. The map $i_1$ is the wedge inclusion in the first factor. Applying $\Sigma^\infty_1$ (the left adjoint to the 1-st space functor from spectra to based spaces) 
gives a homotopy equivalence of $\S$-modules $\sip X \simeq \S \vee \Sigma^\infty X$, and the previous diagram becomes the following commutative diagram of $\S$-modules:
\begin{equation}\label{rew} \xymatrix{& \S \ar[ld]_-e \ar[rd]^-{i_1} \\ \sip X \ar[rr]^-{\simeq} && \S \vee \Sigma^\infty X.}\end{equation}
Here $i_1$ is the wedge inclusion in the first factor.%

Applying this to $X=K(\Z,3)$ and combining it with Theorem \ref{thhku1} and Proposition \ref{kurat}, we obtain weak equivalences of $KU$-modules
\begin{align}\label{goba} THH(KU) &\simeq KU \wedge \sip K(\Z,3) \stackrel{\simeq}{\to} KU \wedge (\S \vee \Sigma^\infty K(\Z,3)) \cong \\ &\cong KU \vee (KU \wedge K(\Z,3)) \simeq KU \vee \Sigma KU_\Q. \nonumber \end{align} %
Note that each of the $KU$-modules in that chain has a map of $\S$-modules from $KU$, namely: $\eta: KU\to THH(KU)$ is the unit, $KU\to KU\wedge \sip K(\Z,3)$ is $\id\wedge e$, $KU\to KU\wedge (\S \vee \Sigma^\infty K(\Z,3))$ is $\id\wedge i_1$, $KU\to KU \vee (KU \wedge K(\Z,3))$ is the inclusion in the first factor and the same goes for $KU \vee \Sigma KU_\Q$. The weak equivalences above are compatible with these maps: the first one because it is a zig-zag of weak equivalences of $KU$-algebras, then we use the commutativity of (\ref{rew}), and then it follows from an inspection of how the distributivity isomorphism works. %

In the homotopy category of $KU\Mod$, we consider the map $\Sigma KU_\Q\to THH(KU)$ which is the inclusion into $KU \vee \Sigma KU_\Q$ followed by the isomorphism obtained from (\ref{goba}). Since $\Sigma KU_\Q$ is a cofibrant $KU$-module, we can represent this map by a morphism of $KU$-modules 
\begin{equation}\label{defff} f:\Sigma KU_\Q\to THH(KU).\end{equation} %

After passing to the homotopy category of $KU$-modules, the map of $KU$-modules \begin{equation}\label{etaf} (\eta,f):KU \vee \Sigma KU_\Q\to THH(KU)\end{equation}
coincides with the isomorphism obtained from (\ref{goba}), by construction and by the remarks right after (\ref{goba}). 
In particular, $(\eta,f)$ is a weak equivalence of $KU$-modules.

The morphism of $KU$-modules $f$ induces a map of commutative $KU$-algebras
\[\tilde f:F(\Sigma KU_\Q) \to THH(KU).\]
To see that it is a weak equivalence, note that, by definition of $\tilde f$, it is such that $\tilde f\circ i=(\eta,f)$, where the weak equivalence of $KU$-modules $i:KU \vee \Sigma KU_\Q\to F(\Sigma KU_\Q)$ was introduced in Remark \ref{remark-i}. 
Since $(\eta,f)$ and $i$ are weak equivalences, then $\tilde f$ is a weak equivalence, too. %
\eprf

\bobs \label{augment-final} In this remark, we prove that $\tilde f:F(\Sigma KU_\Q)\to THH(KU)$ is compatible with the augmentations in a sense to be made explicit below.

Each of the steps in the zig-zag of weak equivalences of $KU$-modules (\ref{goba}) between $KU\wedge K(\Z,3)_+$ and $KU \vee \Sigma KU_\Q$ is augmented over $KU$ and the maps in the zig-zag commute with the augmentations. Here $KU \wedge K(\Z,3)_+$ has the augmentation over $KU$ given by $\id \wedge \epsilon$, as in the last part of Remark \ref{augment-complicado}, and $KU\vee \Sigma KU_\Q$ has the augmentation over $KU$ given by $\id \vee *$, where $*$ denotes the trivial map. As for the intermediate steps, $KU\wedge (\S\vee \Sigma^\infty K(\Z,3))$ is augmented by $\id \wedge (\id \vee *)$, and $KU\vee (KU \wedge K(\Z,3))$ is augmented by $\id \vee (\id \wedge *)=\id\vee *$. The only non-trivial part in this verification is in the first step, the one obtained from the homotopy equivalence $\sip K(\Z,3) \simeq \S \vee \Sigma^\infty K(\Z,3)$ used in (\ref{rew}). More generally, for any well-based space $X$ the following diagram commutes:
\[\xymatrix{\sip X \ar[rr]^-\simeq \ar[rd]_-\epsilon && \S \vee \Sigma^\infty X \ar[ld]^-{\id\vee *} \\ & \S}\]
where $\epsilon:\sip X \to \S$ is obtained from $X\to *$. This follows from the commutativity of the following diagram of based spaces:
\[\xymatrix{\Sigma(X_+) \ar[rr]^-{\simeq} \ar[rd]  && S^1\vee \Sigma X \ar[ld] \\ & S^1}\]
where the left diagonal map is $\Sigma(X\to *)_+$ and the right diagonal map is the identity map on $S^1$ and the suspension of the trivial based map $X\to S^0$ on $\Sigma X$.

Combining this with Remark \ref{augment-complicado}, we obtain a ``diagram''
\[\xymatrix{THH(KU) \ar@{-}[r]^-\sim \ar[d] & KU \wedge \sip K(\Z,3) \ar@{-}[r]^-\sim \ar[d] & KU \vee \Sigma KU_\Q \ar[ld]^-{\id \vee *} \\ KU \ar@{-}[r]^-\sim & KU  }\]
understood to mean that the left square is actually a ladder diagram of commutative $KU$-algebras and the right triangle is a ladder diagram (with the lower side collapsed to $KU$) of $KU$-modules. The triangle commutes (in the sense that all the triangles it hides commute), and the square commutes up to a homotopy of commutative $KU$-algebras. Similarly to how we constructed $f:\Sigma KU_\Q\to KU$, we can obtain a morphism of commutative $KU$-algebras $g:KU\to KU$ such that the following diagram is homotopy commutative in $KU$-modules:
\[\xymatrix{THH(KU) \ar[d] & KU \vee \Sigma KU_\Q \ar[l]_-{(\eta,f)}^-\sim \ar[d]^-{\id\vee *} \\ KU  & KU. \ar[l]^-g_-\sim}\]
Precomposing with the canonical map into the second factor $\Sigma KU_\Q \to KU \vee \Sigma KU_\Q$, we obtain the homotopy commutative diagram in $KU$-modules:
\[\xymatrix{THH(KU) \ar[d] & \Sigma KU_\Q \ar[l]_-{f}^-\sim \ar[d]^-{*} \\ KU  & KU. \ar[l]^-g_-\sim}\]
Since $F$ is left Quillen and $\Sigma KU_\Q$ is a cofibrant $KU$-module, this implies that the diagram
\begin{equation}\label{terf} \xymatrix{THH(KU) \ar[d] & F(\Sigma KU_\Q) \ar[l]_-{\tilde f}^-\sim \ar[d] \\ KU  & KU. \ar[l]^-g_-\sim}\end{equation}
is homotopy commutative in $KU\CAlg$, by an application of \cite[8.5.16]{hirschhorn}.
\eobs

\subsection{The morphism \texorpdfstring{$\sigma$}{sigma}} If $R$ is a commutative $\S$-algebra, there is a natural transformation of $\S$-modules \cite[Section 3]{mc-st}, \cite[IX.3.8]{ekmm}, \cite[3.12]{angeltveit-rognes}\[\sigma: \Sigma R\to THH(R).\]  %

Consider the map of $\S$-modules
\[(\eta,\sigma): KU \vee \Sigma KU \to THH(KU).\] %
It is tempting to conjecture that its rationalization
\[(\eta_\Q,\sigma_\Q):KU_\Q \vee \Sigma KU_\Q \to THH(KU)_\Q\]
is a weak equivalence, since by the results of the previous section, the $\S$-modules $KU_\Q \vee \Sigma KU_\Q$ and $THH(KU)_\Q$ are weakly equivalent. 

However, this is not the case. I thank Geoffroy Horel and Thomas Nikolaus for pointing out this fact and the following proof to me. We will prove that $\sigma:\Sigma KU \to THH(KU)$ is zero in $\pi_1$, therefore it is still zero after rationalization.
By naturality of $\sigma$, we have a commutative diagram
\begin{equation}\label{cuadku}\xymatrix{\Sigma \S \ar[r]^-\sigma \ar[d]_{\Sigma \iota} &THH(\S) \simeq \S \ar[d]^-{THH(\iota)} \\ \Sigma KU \ar[r]_-\sigma & THH(KU)}\end{equation}
where $\iota:\S\to KU$ is the unit of $KU$. After taking $\pi_1$, we obtain a commutative diagram of abelian groups
\begin{equation}\label{cuad}\xymatrix{\Z \ar[r] \ar[d]_-\id & \Z/2 \ar[d] \\ \Z \ar[r] & \Q.}\end{equation}
Therefore, $\Z\to \Q$ must be the zero map, since only the abelian group map $\Z/2\to \Q$ is the zero map.

Note that the same proof works for $L$ (the $p$-adic completion of the Adams summand of $KU$, $p$ a prime) instead of $KU$. Recall that $\pi_*L \cong \Z_{(p)}[(v_1)^{\pm 1}]$, with $v_1$ in degree $2p-2$. %
After replacing $KU$ with $L$ in (\ref{cuadku}) and taking $\pi_1$, we obtain a square which looks like (\ref{cuad}) except with a $\Z_{(p)}$ on the lower left corner. The vertical map $\Z\to \Z_{(p)}$ is the unit of $\Z_{(p)}$: this still forces $\pi_1\sigma: \pi_1(\Sigma L)\to \pi_1(THH(L))$ to be zero. \index{L}

This corrects an error in \cite[8.4]{mc-st} where it is claimed that there is a weak equivalence $L_\Q \vee \Sigma L_\Q \stackrel{\sim}{\to} THH(L)_\Q$ induced by $(\eta,\sigma)$. As a positive result, we have Theorem \ref{thhku2} and the weak equivalence (\ref{etaf}) instead.

\section{Iterated topological Hochschild homology of \texorpdfstring{$KU$}{KU}} \label{sect-iterated} 

Let $A$ be a commutative $\S$-algebra. We denote by $THH^n(A)$ the \emph{iterated topological Hochschild homology of $A$}, i.e. $THH(\dots(THH(A)))$ where $THH$ is applied $n$ times. Other expressions for $THH^n(A)$ include $T^n \otimes A$ or $\Lambda_{T^n}(A)$, where $T^n$ is an $n$-torus and $\Lambda$ is the Loday functor \cite{cdd}. Note that $THH^n(A)$ is an augmented commutative $A$-algebra.

We will now give two different descriptions of $THH^n(KU)$ for $n\geq 2$. The first one, given in Theorem \ref{thhnku1cor}, generalizes Theorem \ref{thhku1} which describes $THH(KU)$ via Eilenberg-Mac Lane spaces. The second one, given in Theorem \ref{thhnkufree}, generalizes Theorem \ref{thhku2} which describes $THH(KU)$ as a free commutative $KU$-algebra on a $KU$-module.

We have also given a description of the commutative $KU$-algebra $THH(KU)$ as the split square-zero extension $KU \vee \Sigma KU_\Q$ in Theorem \ref{thhku2}. For $n\geq 2$, $THH^n(KU)$ is not a split square-zero extension of $KU$, as we shall see. However, it is a split extension: we will describe the non-unital commutative algebra structure of the homotopy groups of its augmentation ideal, which is rational as in the $n=1$ case.

\subsection{Description via Eilenberg-Mac Lane spaces} Let $G$ be a CW-complex which is a topological abelian group with unit a 0-cell and a cellular addition map. Applying Lemma \ref{sinloc} and Proposition \ref{Smono}, we obtain isomorphisms of commutative $\S[G]$-algebras: 
\begin{align*}THH^2(\S[G])& \cong THH(\S[G] \wedge \S[BG]) \cong THH(\S[G\times BG]) \\ & \cong \S[G\times BG] \wedge \S[B(G\times BG)] \cong \S[G] \wedge \S[BG \times BG \times B^2G] \end{align*}
which we have written as $\S[G][BG \times BG \times B^2G]$.

For general $n\geq 2$, the same type of computation gives a description of $THH^n(\S[G])$: we obtain an isomorphism of commutative $\S[G]$-algebras
\begin{equation}\label{ane}THH^n(\S[G]) \cong \S[G][B^{a_1}G\times \dots\times B^{a_{2^n-1}}G].\end{equation} The numbers $a_i$ can be described as follows. Let $v_0=0$. Define by induction \begin{equation}\label{defai} v_n=(v_{n-1},v_{n-1}+(1,\dots,1)) = (a_0, \dots,a_{2^n-1})\in \N^{2^n}\end{equation}
for $n\geq 1$. For example, $v_1=(0,1)$, $v_2=(0,1,1,2)$ and $v_3=(0,1,1,2,1,2,2,3)$. %
This sequence of integers can be found in the On-Line Encyclopedia of Integer Sequences \cite{oeis}. We can give an easier description. Let $I_n$ be the multiset having as elements the numbers $i$ with multiplicity ${n\choose i}$, for $i=1,\dots,n$. Denote the multiplicity of an element $x$ of a multiset by $|x|$. Now note that the multiset underlying the sequence $(a_1,\dots,a_{2^n-1})$ defined in $(\ref{defai})$ coincides with $I_n$, by Pascal's rule. %
Therefore, the isomorphism (\ref{ane}) can be reformulated as
\begin{equation}\label{thhnsg-ref} THH^n(\S[G])\cong \S[G]\left[\prod\limits_{i=1}^n (B^iG)^{\times {n \choose i}}\right].\end{equation}

The following theorem generalizes Theorem \ref{teoloc} to higher iterations of $THH$.

\bteo \label{thhnloc} Let $x\in \pi_*\S[G]$ and $n\geq 1$. There is a zig-zag of weak equivalences of commutative $\S[G][x^{-1}]_h$-algebras
\begin{equation}\label{gtre}THH^n(\S[G][x^{-1}]_h) \simeq \S[G][x^{-1}]_h[B^{a_1}G\times \dots\times B^{a_{2^n-1}}G],\end{equation}
or alternatively,
\begin{equation}\label{gtre2}THH^n(\S[G][x^{-1}]_h)\simeq  \S[G][x^{-1}]_h\left[\prod\limits_{i=1}^n (B^iG)^{\times {n \choose i}}\right].\end{equation}
\bprf The proof is by induction. The base case is Theorem \ref{teoloc}. We do the induction step for $n=2$ for simplicity: for higher $n$ it is analogous, only more cumbersome to write down.

We use the notations introduced for the proof of Theorem \ref{teoloc}, namely $A=\S[G]$, $B=\S[BG]$, $Q$ is the cofibrant replacement functor in $\SCAlg$ obtained from a given functorial factorization in $\SCAlg$, and $Q_A$ is the cofibrant replacement functor in $QA\CAlg$ obtained by factoring the unit of a $QA$-commutative algebra in $\SCAlg$. From the proof of that theorem, we obtain a zig-zag of weak equivalences of cofibrant commutative $QA[x^{-1}]$-algebras which we follow by an isomorphism gotten by applying Proposition \ref{Smono}:
\begin{equation}\label{gohj} THH(QA[x^{-1}]) \simeq Q_A(A\wedge B)[(x\wedge 1)^{-1}] \cong Q_A \S[G \times BG][(x\times 1)^{-1}].\end{equation}
We apply $THH$ to obtain weak equivalences of commutative $QA[x^{-1}]$-algebras
\begin{align*}THH^2(QA[x^{-1}]) &\simeq THH(Q_A\S[G\times BG][(x\times 1)^{-1}]) \\
&\simeq (Q_A \S[G \times BG])[(x\times 1)^{-1}] \wedge Q\S[B(G\times BG)].\end{align*}
The second line comes from the zig-zag of weak equivalences obtained as follows. First, we observe that if $C$ is a $QA$-commutative algebra, then $Q_AC\to C$ defines a functorial cofibrant replacement of $C$ in the category of commutative $\S$-algebras. %
We apply this remark to $C=\S[G\times BG]$. As a consequence, analogous steps to the ones taken in the proof of Theorem \ref{teoloc} from the beginning up to (\ref{vers}) apply \emph{mutatis mutandis} and get us the result. 
We continue:
\begin{align*}
(Q_A \S[G \times BG])[(x\times 1)^{-1}] \wedge Q\S[B(G\times BG)] &\cong Q_A(A\wedge B)[(x\wedge 1)^{-1}] \wedge Q\S[BG\times B^2G] \\
&\simeq QA[x^{-1}]\wedge QB \wedge Q\S[BG\times B^2G]\\
&\simeq (QA[x^{-1}])[BG\times BG \times B^2G].
\end{align*}
Here, the first step is gotten from the isomorphism in (\ref{gohj}). In the course of the proof of Theorem \ref{teoloc} we obtained a zig-zag of weak equivalences of commutative $QA[x^{-1}]$-algebras $Q_A(A\wedge B)[(x\wedge 1)^{-1}] \simeq QA[x^{-1}]\wedge QB$. The steps in this zig-zag are all cofibrant commutative $\S$-algebras and $Q\S[BG\times B^2G]$ is a cofibrant commutative $\S$-algebra: this explains the second step. The third step is an application of Lemma \ref{cofcofsemi} together with the isomorphism from Proposition \ref{Smono}.
\eprf
\eteo

\bobs \label{augment-complicado2} As in Remark \ref{augment-complicado}, the zig-zag in the previous theorem is compatible up to homotopy of commutative $QA[x^{-1}]$-algebras with the augmentations, where the right-hand side of (\ref{gtre}) is augmented by means of the map $B^{a_1}G\times \dots\times B^{a_{2^n-1}}G\to *$ and similarly in (\ref{gtre2}).
\eobs

As a corollary, we obtain:

\bteo \label{thhnku1cor}
There is a zig-zag of weak equivalences of commutative $KU$-algebras
\begin{equation}\label{thhnku1}THH^n(KU)\simeq KU[K(\Z,a_1+2) \times \dots \times K(\Z,a_{2^n-1}+2)],\end{equation}
or alternatively,
\begin{equation}\label{thhnku2}THH^n(KU)\simeq KU\left[\prod\limits_{i=1}^n K(\Z,i+2)^{\times {n \choose i}}\right].\end{equation}
\eteo
For example,
\begin{equation}\label{thh2kua} THH^2(KU)\simeq KU[K(\Z,3)\times K(\Z,3)\times K(\Z,4)].\end{equation}

The previous theorem generalizes the expression of Theorem \ref{thhku1} for $THH(KU)$ as \linebreak $KU[K(\Z,3)]$ to $THH^n(KU)$. We can also generalize the expression for $THH(KU)$ as $F(\Sigma KU_\Q)$ of Theorem \ref{thhku2}. Note that the proof uses results from Section \ref{section:snku} below.

\bteo \label{thhnkufree} Let $n\geq 2$. There is a zig-zag of weak equivalences of commutative $KU$-algebras
\[F\left(\bigvee\limits_{i=1}^n (S^i)^{\vee {n\choose i}} \wedge KU_\Q\right)
\simeq THH^n(KU).\]
\bprf Since $- \wedge KU_\Q:\SMod\to KU\Mod$ and $F:KU\Mod\to KU\CAlg$ are left adjoints, they preserve coproducts, so:
\[F\left(\bigvee\limits_{i=1}^n (S^i)^{\vee {n\choose i}} \wedge KU_\Q\right)
\cong F\left(\bigvee_{i=1}^n (\Sigma^i KU_\Q)^{\vee {n\choose i} } \right) \cong \bigwedge_{\substack{KU \\ i=1}}^n F(\Sigma^i KU_\Q)^{\wedge_{KU} {n\choose i}}.\]
From (\ref{snfree}), we obtain a zig-zag of weak equivalences of commutative $KU$-algebras %
\[\bigwedge_{\substack{KU \\ i=1}}^n F(\Sigma^i KU_\Q)^{\wedge_{KU} {n \choose i}} \simeq  \bigwedge_{\substack{KU \\ i=1}}^n (S^i\otimes KU)^{\wedge_{KU} {n \choose i}}. \] %
Using Theorem \ref{holo}, %
\[\bigwedge_{\substack{KU \\ i=1}}^n (S^i\otimes KU)^{\wedge_{KU} {n \choose i}} \simeq \bigwedge_{\substack{KU \\ i=1}}^n KU[K(\Z,i+2)]^{\wedge_{KU} {n \choose i}} \cong KU\left[ \prod_{i=1}^n K(\Z,i+2)^{\times {n\choose i}}\right]\]
which is weakly equivalent to $THH^n(KU)$ by Theorem \ref{thhnku1cor}.
\eprf
\eteo

\bobs \label{stabeq} We might be tempted to prove the previous theorem more directly, arguing from the weak equivalence of based spaces
\begin{equation}\label{stabletorus}\Sigma T^n\simeq \Sigma \bigvee\limits_{i=1}^n (S^i)^{\vee {n \choose i}}.\end{equation} However, we do not know a priori whether this guarantees that $THH^n(KU)= T^n\otimes KU$ is weakly equivalent to $\left( \bigvee\limits_{i=1}^n (S^i)^{\vee {n \choose i}}\right) \otimes KU$ (which we can easily compute using the description from Theorem \ref{holo} of $S^i \otimes KU$ for all $i\geq 1$ and the fact that $-\otimes KU$ preserves coproducts). Indeed, there are counterexamples to the statement that if $A$ is a commutative $\S$-algebra, then $X\otimes A \simeq Y\otimes A$ provided $\Sigma X \simeq \Sigma Y$ \cite{dundas-tenti}. After having proved the theorem, though, we have that $KU$ does satisfy this for the special case of (\ref{stabletorus}). The partial results of \cite{veen} (extended by \cite{blprz}) prove that $A=H\F_p$ also satisfy it for (\ref{stabletorus}), at least in a certain range relating $n$ and $p$. 
We are led to ask ourselves the question, as \cite[4.1]{dundas-tenti} did for $A=H\F_p$, of whether more generally $KU$ is such that $X\otimes KU \simeq Y\otimes KU$ provided $\Sigma X\simeq \Sigma Y$. More ambitiously, it would be interesting to find conditions on any commutative $\S$-algebra $A$ that guarantee this property.\eobs

\subsection{The augmentation ideal} \label{second-desc}In this section, we investigate the augmentation ideal of \linebreak $THH^n(KU)$. Let us define this concept.

Let $R$ be a commutative $\S$-algebra and $A$ be a commutative $R$-algebra with augmentation $\epsilon:A\to R$. Denote by $\overline A$ the fiber of $\epsilon$, i.e. it is the $R$-module obtained as the pullback
\[\xymatrix{\overline A \ar[d] \ar[r]^-i & A \ar[d]^-\epsilon \\ {*} \ar[r] & R}\] 
in $\RMod$. It gets a non-unital multiplication from the universal property of pullbacks, by considering the following commutative diagram in $\RMod$. See \cite[Section 2]{basterra} for further elaboration.
\[\xymatrix{\overline A \wedge_R \overline A \ar[d] \ar[r]^{i\wedge i} & A \wedge_R A \ar[d]^-{\epsilon \wedge \epsilon} \ar[r]^-\mu & A \ar[d]^-\epsilon \\ {*} \ar[r] & R \wedge_R R \ar[r]_-\cong & R}\]

Consider the augmentation $\epsilon:THH^n(KU)\to KU$. To ensure that we compute its \emph{homotopy} fiber, we first replace $\epsilon$ by a fibration in the category of commutative $KU$-algebras, i.e. we replace $\epsilon$ by the fibration appearing in its factorization by an acyclic cofibration followed by a fibration. We denote the fiber of this new fibration by $\overline{THH}^n(KU)$.

We first need a generalization of Proposition \ref{kurat}:

\begin{prop} \label{lemgen} Let $r\geq 3$. There are zig-zags of weak equivalences of $KU$-modules
\[KU \wedge K(\Z,r) \simeq \begin{cases} \Sigma KU_\Q & \text{if } r \text{ is odd}, \\ \bigvee\limits_{m\geq 1} KU_\Q & \text{if } r \text{ is even}.\end{cases}\]
\bprf When $r$ is odd, the proof of Proposition \ref{kurat} works just as well, and when $r$ is even it gives us
\[KU \wedge K(\Z,r) \simeq KU_\Q \wedge K(\Z,r)_\Q.\]
So let $r$ be even. As noted in Section \ref{rationalization}, $K(\Z,r)_\Q \simeq \Omega S^{r+1}_\Q$. Now we use the James splitting which says that, for $X$ a connected based CW-complex, $\Sigma \Omega \Sigma X \simeq \Sigma \bigvee_{m\geq 1} X^{\wedge m}$. %
Therefore, %
$\Sigma^\infty \Omega \Sigma X \simeq \Sigma^\infty \bigvee_{m\geq 1} X^{\wedge m}$. Rationalizing it and applying it to $X=S^r$, we obtain
\[\Sigma^\infty K(\Z,r)_\Q \simeq \Sigma^\infty \Omega S^{r+1}_\Q \simeq \Sigma^\infty \bigvee_{m\geq 1} S^{rm}_\Q.\] 
Since $r$ is even, Bott periodicity gives the result.
\eprf
\end{prop}

\bcor \label{hocs} The augmentation ideal $\overline{THH}^n(KU)$ is rational.
\bprf The expression (\ref{thhnku1}) gives, after splitting off the units of the spherical group rings like in (\ref{goba}), a zig-zag of weak equivalences of $KU$-modules
\begin{equation} \label{kol} THH^n(KU) \simeq KU \wedge (\S \vee \Sigma^\infty K(\Z,a_1+2)) \wedge \dots \wedge  (\S \vee \Sigma^\infty K(\Z,a_{2^n-1}+2)). \end{equation}
Distributing the terms in the right-hand side and applying Proposition \ref{lemgen} proves that the homotopy fiber of the augmentation of the right-hand side is rational. 
Since the zig-zag is compatible up to homotopy with the augmentations (Remark \ref{augment-complicado2}), this implies that the homotopy fibers are weakly equivalent, so in particular $\overline{THH}^n(KU)$ is rational.
\eprf
\ecor

Recall that $H\Q$ is a homotopy commutative ring spectrum whose multiplication map is a weak equivalence, and it is a cofibrant $\S$-module. In particular, $KU_\Q$ is a homotopy  commutative ring spectrum, and from (\ref{thhnku1}) we obtain a zig-zag of weak equivalences of commutative $KU_\Q$-ring spectra:
\[THH^n(KU)_\Q \simeq KU_\Q \wedge K(\Q,a_1+2)_+ \wedge \dots \wedge K(\Q,a_{2^n-1}+2)_+.\]
By using the identification of the rationalized Eilenberg-Mac Lane spaces of Section \ref{rationalization} and the computation of the rational homology of loop spaces of odd-dimensional spheres \cite[Page 225]{fht-rational}, we obtain
\bprop There is an isomorphism of commutative $\Q[t^{\pm 1}]$-algebras
\begin{equation} \label{thhnkuq} H\Q_*(THH^n(KU)) \cong \Q[t^{\pm1}] \otimes \bigotimes_{a_i \textup{ odd}} E(\sigma^it) \otimes \bigotimes_{a_j \textup{ even}} \Q[\sigma^j t]  \end{equation} %
where $|\sigma^rt|=a_r+2$ and $i,j \in\{1,\dots, 2^n-1\}$. \eprop
For example,
\[H\Q_*(THH^2(KU)) \cong \Q[t^{\pm 1}] \otimes E(\sigma t) \otimes E(\sigma t)  \otimes \Q[\sigma^2t]\]
with $|\sigma t|=3$ and $|\sigma^2t|=4$.\\ %

We can recognize the right-hand side of the expression (\ref{thhnkuq}) as an iterated Hochschild homology algebra: \index{Hochschild homology!iterated}
\begin{equation} H\Q_*(THH^n(KU)) \cong HH^{\Q,n}_*(\Q[t^{\pm 1}]).\end{equation}
Indeed, $HH^\Q_*(\Q[t^{\pm1}])\cong \Q[t^{\pm 1}] \otimes E(\sigma t)$, and $HH^\Q_*(E(\sigma t))\cong E(\sigma t)\otimes \Q[\sigma^2t]$. These Hochschild homology calculations are classical and can be found e.g. in \cite[Section 2]{mc-st} and \cite[2.4]{angeltveit-rognes}. We use that localization commutes with Hochschild homology \cite[Theorem 9.1.8(3)]{weibel}. Also note that in general, the Hochschild homology of an exterior algebra is isomorphic to the tensor product of this same exterior algebra with a divided power algebra, but over $\Q$ such algebras are polynomial.\\

Denote by $\overline{HH}_*^{\Q,n}(B)$ the kernel of the augmentation $HH_*^{\Q,n}(B)\to B$. From these remarks and the proof of Corollary \ref{hocs}, we obtain:
\bteo  \label{high-hoch} There is an isomorphism of non-unital commutative $\Q[t^{\pm 1}]$-algebras
\[\overline{THH}^n_*(KU) \cong \overline{HH}_*^{\Q,n}(\Q[t^{\pm 1}]).\]
\eteo

\section{\texorpdfstring{$\Sigma Y\otimes KU$}{Sigma Y tensor KU}} \label{section:snku}

In this section, we evaluate the commutative $KU$-algebra $\Sigma Y \otimes KU$ when $Y$ is a based CW-complex, by comparing it with $Y\otimes_{KU}(S^1\otimes KU)$. We are very grateful to Bj{\o}rn Dundas for suggesting this line of argument.\\ 

Recall that if $R$ is a commutative $\S$-algebra, the category $\RCAlg$ is tensored over $\Top$ \cite[VII.2.9]{ekmm}. If $A\in \RCAlg$, then the tensor $S^1\otimes_R A$ is naturally isomorphic to $THH^R(A)$ as a commutative augmented $A$-algebra \cite{mc-schw-vo}, \cite[IX.3.3]{ekmm}, \cite[Section 3]{angeltveit-rognes}. Therefore, in this section we will identify $S^1\otimes_R A$ and $THH^R(A)$ without further notice.

\subsection{The morphism \texorpdfstring{$\nu$}{nu}}
\def\C{\mathcal{C}}
Let $\C$ be a category enriched and tensored over $\Top$. Denote its tensor by $\otimes$. Fix a based space $(Z,z_0)$. We denote by $\nu^Z$ the natural transformation
\begin{equation}\label{nuu}\xymatrix@C+2pc{\C \rtwocell<5>^{\id}_{Z\otimes -}{\;\;\nu^Z}  & \C}\end{equation}
whose component in $C\in \C$ is given by \begin{equation}\label{defnu}\nu^Z_C\coloneqq\eta_Z^C(z_0):C\to Z\otimes C.\end{equation} Here $\eta_Z^C:Z\to \C(C,Z\otimes C)$ is the unit at $Z$ of the adjunction 
\begin{equation} \label{adjo} \xymatrix@C+2pc{\Top \ar@/^.8pc/[r]^(.5){-\otimes C} & \C. \ar@/^.8pc/[l]^-{\C(C,-)}}
\end{equation}

Let us now highlight the naturality properties of $\nu^Z_C$ at $C$ and at $Z$. Let $\varphi:C\to C'$ be a morphism in $\C$. The naturality of the isomorphism \[\C(Z\otimes C,Z\otimes -)\cong \Top(Z,\C(C,Z\otimes -))\] gives the commutativity of the following diagram.
\begin{equation}\label{cuadradonu}\xymatrix{C\ar[d]_-\varphi \ar[r]^-{\nu^Z_C} & Z\otimes C \ar[d]^-{\id\otimes \varphi} \\ C' \ar[r]_-{\nu_{C'}^Z} & Z\otimes C'}\end{equation}
Let $u:Z\to Z'$ be a morphism of based spaces. The naturality of $\eta^C$ gives the commutativity of the following diagram.
\begin{equation}\label{nuu2}\xymatrix{C\ar[r]^-{\nu_C^Z} \ar[rd]_-{\nu_C^{Z'}} & Z\otimes C \ar[d]^-{u\otimes \id} \\ & Z'\otimes C }\end{equation}

\bej The category $\Top_*$ is tensored over $\Top$: if $X\in \Top_*$ and $Y\in \Top$, then $Y\otimes X$ is defined as $Y_+\wedge X$. When $(Y,y_0)$ is based, we denote by \begin{equation}\label{nyx} n^Y_X:X\to Y_+\wedge X\end{equation} the map $\nu_X^Y$ of (\ref{defnu}) applied to $\C=\Top_*$. More explicitely, the map $n^Y_X$ takes $X$ to the copy of $X$ lying over $y_0$ in $Y_+\wedge X$. %
\eej

\subsection{In commutative algebras}
Let $R$ be a commutative $\S$-algebra. Let $A$ be a commutative $R$-algebra and $(X,x_0)$ be a based space. The map (\ref{defnu}) in this scenario is a map of commutative $R$-algebras \[\nu^X_A:A\to X\otimes_R A\] which gives $X\otimes_R A$ the structure of a commutative $A$-algebra. In particular, when $X=S^1$, this is the usual structure of an $A$-algebra of $THH^R(A)$.%

Now, take $R=\S$ and $A=KU$. Let $(Y,y_0)$ be a based space. We use the symbol $\otimes$ to denote the tensor of $\SCAlg$ over $\Top$. Consider the following diagram in $KU\CAlg$. Here the map $e:S^1\to *$ collapses the circle into its basepoint, and we have identified $F(*\wedge KU_\Q)$ and $*\otimes KU$ with $KU$. 
\begin{equation}\label{diagorra}\xymatrix@C+2pc{
Y\otimes_{KU} F(S^1\wedge KU_\Q) \ar[d]^-{\id \otimes \tilde f}_\sim & F(S^1\wedge KU_\Q) \ar[l]_-{\nu^Y_{F(S^1\wedge KU_\Q)}} \ar[d]^-{\tilde f}_\sim \ar[r]^-{F(e\wedge \id)} & KU \ar[d]^-g_-\sim \\
Y\otimes_{KU} (S^1\otimes KU)  & S^1\otimes KU  \ar[l]^-{\nu^Y_{S^1\otimes KU}} \ar[r]_{e\otimes \id} & KU }\end{equation}
The weak equivalence $\tilde f$ comes from Theorem \ref{thhku2}. The map $g:KU\to KU$ comes from (\ref{terf}): in that remark we proved that the right square commutes up to a homotopy of commutative $KU$-algebras. The left square commutes as an application of the commutativity of (\ref{cuadradonu}). Note that $\id\otimes \tilde f$ is a weak equivalence because $Y\otimes_{KU}-$ is a left Quillen functor, assuming $Y$ is a based CW-complex.

We will now identify the members of the left column. 

\begin{prop} Let $(X,x_0)$ and $(Y,y_0)$ be based spaces, and let $A$ be a commutative $R$-algebra. %
\be
\item There is an isomorphism of commutative $A$-algebras
\[Y \otimes_A (X \otimes_R A) \cong (Y_+ \wedge X) \otimes_R A\]
where $\otimes_R$ (resp. $\otimes_A$) denotes the tensoring of $\RCAlg$ (resp. $\ACAlg$) over $\Top$.%

Moreover, the isomorphism makes the following diagram in $A\CAlg$ commute. The morphism $n^Y_X:X\to Y_+\wedge X$ was defined in (\ref{nyx}).
\begin{equation}\label{ndiag}\xymatrix@C+1pc{X\otimes_R A \ar[r]^-{\nu^Y_{X\otimes_RA}} \ar[rd]_-{n^Y_X\otimes \id}  & Y\otimes_A (X\otimes_R A) \ar[d]^-\cong \\ & (Y_+\wedge X) \otimes_R A}\end{equation}
\item Let $M$ be an $A$-module. Let $F:A\Mod\to A\CAlg$ be the free commutative algebra functor. There is an isomorphism
\[Y\otimes_A F(X\wedge M) \cong F(Y_+ \wedge X\wedge M)\]
making the following diagram commute.
\[\xymatrix@C+1pc{F(X\wedge M) \ar[rd]_-{F(n^Y_X \wedge \id)} \ar[r]^-{\nu^Y_{F(X\wedge M)}} & Y\otimes_A F(X\wedge M) \ar[d]^-\cong \\ & F(Y_+\wedge X\wedge M)}\] 
In the expression $Z\wedge M$ for a based space $Z$ we are using the tensor of $\AMod$ over $\Top_*$.%
\ee
\bprf (1) Let $B$ be a commutative $A$-algebra with unit $\varphi: A\to B$. Using the defining adjunction for $Y\otimes_A-$, we get a homeomorphism
\begin{equation}\label{pro1}A\CAlg(Y\otimes_A (X \otimes_R A),B) \cong \Top(Y,A\CAlg(X\otimes_R A,B)).\end{equation}
The morphisms of commutative $A$-algebras $X\otimes_R A\to B$ are the morphisms of commutative $R$-algebras $g:X\otimes_R A\to B$ making the following diagram commute:
\[\xymatrix{& A \ar[rd]^-\varphi \ar[ld]_-{\nu_A^{X}} \\ X\otimes_R A \ar[rr]_-g  && B.}\]
Recalling the definition of $\nu$, this means that 
\begin{equation}\label{hoji}g\circ \eta^A_{X}(x_0)=\varphi.\end{equation}
The adjoint map of $g$ by the defining adjunction of $-\otimes_R A$ is the map in $\Top$
\begin{equation}\label{hojiadj}\xymatrix{X\ar[r]^-{\eta_{X}^A} & \RCAlg(A,X\otimes_R A) \ar[r]^-{g_*} & \RCAlg(A,B).}\end{equation}
Let the space $\RCAlg(A,B)$ be pointed by $\varphi:A\to B$. The condition (\ref{hoji}) on the map $g$ is then translated to the adjoint (\ref{hojiadj}) by stating that it is a based map, i.e. it takes $x_0$ to $\varphi$. Thus, continuing (\ref{pro1}),
\begin{equation}\label{prox}\Top(Y,\ACAlg(X\otimes_R A,B)) \cong \Top(Y,U \Top_*(X,\RCAlg(A,B))),\end{equation}
where $U:\Top_*\to \Top$ is the functor forgetting the basepoint. It is the right adjoint to the functor $(-)_+:\Top\to \Top_*$ which adds a disjoint basepoint, so we continue:
\[\Top(Y,U \Top_*(X,\RCAlg(A,B))) \cong U\Top_*(Y_+,\Top_*(X,\RCAlg(A,B))).\] %
Since $\Top_*(X,-):\Top_*\to \Top_*$ is the right adjoint to $-\wedge X$, we get:
\[U\Top_*(Y_+,\Top_*(X,\RCAlg(A,B))) \cong U\Top_*(Y_+\wedge X,\RCAlg(A,B)).\]
By the same argument proving (\ref{prox}), we get
\[U\Top_*(Y_+\wedge X,\RCAlg(A,B)) \cong \ACAlg((Y_+ \wedge X) \otimes_R A,B).\]
In conclusion, we have a homeomorphism
\[A\CAlg(Y\otimes_A (X \otimes_R A),B) \cong \ACAlg((Y_+ \wedge X) \otimes_R A,B),\]
and the Yoneda lemma finishes the proof.

The isomorphism was established using a chain of adjunctions. Following this chain, one observes that both $n^Y_X$ and $\nu_{X\otimes_RA}^Y$, which are defined via units of adjunctions by analogous procedures, make the diagram (\ref{ndiag}) commute.\\

(2) The functor $F$ is defined via a \emph{continuous} monad in $\AMod$ (i.e. it is enriched over $\Top$), see \cite[proof of VII.2.9]{ekmm}. Therefore, the functor $F$ preserves tensors over $\Top$, so we get the desired isomorphism. %
\eprf
\end{prop}

Applying the previous proposition to $R=\S$, $A=KU$, $X=S^1$ and $M=KU_\Q$, the diagram (\ref{diagorra}) can be replaced with the following one.
\begin{equation}\label{diagorra2}\xymatrix@C+2pc{
F(Y_+ \wedge S^1\wedge KU_\Q) \ar[d]^-\sim & F(S^1 \wedge KU_\Q) \ar[l]_-{F(n^Y_{S^1}\wedge \id)} \ar[r]^-{F(e\wedge \id)} \ar[d]^-{\tilde f}_\sim & KU \ar[d]^-g_-\sim \\
(Y_+\wedge S^1) \otimes KU   & S^1\otimes KU \ar[l]^-{n^Y_{S^1}\otimes \id} \ar[r]_-{e\otimes \id}  & KU  
 }\end{equation}
We suppose that $Y$ is a based CW-complex, so that the vertical map on the left is a weak equivalence.

Now, note that the following is a pushout square of based or unbased spaces. %
\[\xymatrix{S^1\ar[r]^-e \ar[d]_-{n^Y_{S^1}} & \ast \ar[d] \\ Y_+\wedge S^1 \ar[r] & Y\wedge S^1}\]
Since the functors $-\otimes KU:\Top\to KU\CAlg$ and $F(-\wedge KU_\Q):\Top_*\to KU\CAlg$ are left adjoints, they preserve pushouts, so the pushout of the top line of (\ref{diagorra2}) is $F(Y\wedge S^1\wedge KU_\Q)$ and the pushout of the bottom line is $(Y\wedge S^1) \otimes KU.$
Now, the three vertical maps of (\ref{diagorra2}) are weak equivalences. The horizontal maps pointing left are cofibrations: indeed, $n^Y_{S^1}$ is a cofibration, %
 $KU_\Q$ is a cofibrant $KU$-module (similarly as in Remark \ref{f-cof}) so $- \wedge KU_\Q$ is left Quillen, %
  $F$ is left Quillen and $-\otimes KU$ is left Quillen. Moreover, all the objects in the diagram are cofibrant in $KU\CAlg$. Since the left square commutes and the right square commutes up to a homotopy of commutative $KU$-algebras, an application of Lemma \ref{hompus} plus the naturality of $n^Y_{S^1}$ in $Y$ (\ref{nuu2}) proves the following
\bteo \label{xku}There is a zig-zag of weak equivalences of commutative $KU$-algebras \[F(Y\wedge S^1\wedge KU_\Q) \simeq (Y\wedge S^1)\otimes KU\]
natural in the based CW-complex $Y$.
\eteo

This determines $\Sigma Y\otimes KU$ as the free commutative $KU$-algebra on the $KU$-module $\Sigma Y \wedge KU_\Q$, up to weak equivalence. In particular, we have a zig-zag of weak equivalences of commutative $KU$-algebras
\begin{equation} \label{snfree} F(\Sigma^n KU_\Q) \simeq S^n \otimes KU\end{equation}
for every $n\geq 1$. %

As in Remark \ref{f-cof}, the $KU$-modules $\Sigma^n KU_\Q$ are cofibrant for $n\geq 0$. Since $F$ is a left Quillen functor, Bott periodicity implies that we have zig-zags of weak equivalences of commutative $KU$-algebras
\begin{equation}S^n\otimes KU \simeq \begin{cases} F(\Sigma KU_\Q) & \text{if } n \text{ is odd,} \\ F(KU_\Q) & \text{if } n \text{ is even} \end{cases}\end{equation}
for every $n\geq 1$.

The line (\ref{snfree}) generalizes the expression of Theorem \ref{thhku2} for $THH(KU)$ as the free commutative $KU$-algebra on $\Sigma KU_\Q$. The following generalizes the expression of Theorem \ref{thhku1} for $THH(KU)$ via Eilenberg-Mac Lane spaces.

\bteo \label{holo}Let $n\geq 1$. Then $S^n \otimes KU \simeq KU[K(\Z,n+2)]$ as commutative $KU$-algebras. %
\eteo
\bprf %
We learned of results similar to the following from \cite{veen}: if $\S\to A\to B$ are cofibrations of commutative $\S$-algebras, then there is a weak equivalence \[S^{n+1}\otimes_A B \stackrel{\sim}{\leftarrow} B^A(B,S^n\otimes_A B,B)\]
where the term on the right side is a two-sided bar construction. Here $\otimes_A$ denotes the tensor of commutative $A$-algebras over $\Top$. Let us give a proof. Since the functor $-\otimes_A B$ is left Quillen, it preserves pushouts and cofibrations, so we have a pushout of commutative $A$-algebras where the arrows $S^n\otimes_A B \to D^{n+1}\otimes_A B$ are cofibrations:
\[\xymatrix{S^n\otimes_A B \ar@{>->}[r] \ar@{>->}[d] & D^{n+1}\otimes_A B \ar[d] \\ D^{n+1}\otimes_A B \ar[r] & S^{n+1}\otimes_A B.}\]
Therefore, 
\begin{align*}
S^{n+1}\otimes_A B & \cong (D^{n+1}\otimes_A B) \wedge_{S^n\otimes_A B} (D^{n+1}\otimes_A B) \\ &\stackrel{\sim}{\leftarrow} B^A(D^{n+1}\otimes_A B,S^n \otimes_A B,D^{n+1}\otimes_A B)\\
&\stackrel{\sim}{\leftarrow} B^A(B, S^n\otimes_A B, B)\end{align*}
where the weak equivalence in the middle is an application of \cite[VII.7.3]{ekmm}, %
and the last one comes from two applications of \cite[VII.7.2]{ekmm}. %

We use this to prove the result by induction. The result is true for $n=1$ (Theorem \ref{thhku1}); suppose it is true for some $n\geq1$. Then
\begin{align*}
S^{n+1}\otimes KU &\simeq B^\S(KU,KU[K(\Z,n+2)],KU) \\
&\simeq B^\S(KU,KU,KU) \wedge B^\S(\S,\S[K(\Z,n+2)],\S) \\ %
&\simeq KU \wedge \S[K(\Z,n+3)] = KU[K(\Z,n+3)]. %
\end{align*}
Here we have used that $B^\S(\S,\S[G],\S)\cong \S[BG]$ for $G$ a topological commutative monoid. This result is proven in the same fashion as Proposition \ref{conmutarTHH}, which deals with the analogous result for the \emph{cyclic} bar construction.
\eprf

\bobs \label{x-thom} In Remark \ref{thh-thom} we observed that $KU$ behaves like a Thom spectrum to the eyes of topological Hochschild homology. Comparing Theorems \ref{thhnku1cor} and \ref{holo} with \cite[1.1]{schlichtkrull-higher} or \cite[4.11]{rsv-thom}, we see that, more generally, $KU$ behaves like a Thom spectrum to the eyes of $X\otimes -$ when $X$ is an $n$-torus or an $n$-sphere, $n\geq 1$. See also Remark \ref{jurs}.\ref{taq-ku-thom} for a similar observation about $TAQ$.\eobs

\section{Topological André-Quillen homology of \texorpdfstring{$KU$}{KU}} \label{section:taq}

If $A\to B$ is a morphism of commutative $\S$-algebras, one can define its \emph{cotangent complex} $\Omega_{B|A}\in B\Mod$, also known as its \emph{topological André-Quillen} $B$-module, $TAQ(B|A)$: see \cite{basterra}. %
We adopt the latter notation. When $A=\S$, we delete it from the notation.

\bteo \label{taqku1} The $KU$-modules $TAQ(KU)$ and $KU \wedge \K$ are weakly equivalent.
\eteo
Here $\K$ is the $\S$-module associated to the topological abelian group $K(\Z,2)$: it is a model for $\Sigma^2 H \Z$. %
More generally, as explained in \cite{basterra-mandell} before Theorem 5, for a topological abelian group $G$ there is an $\S$-module associated to $G$ whose zeroth space is $G$. We denote it by $\mathbf{G}$. More generally, we denote by $\mathbf{X}$ the $\S$-module associated to an $E_\infty$-space $X$ whose zeroth space is the group completion of $X$.

In the next proof we will use the localization of a module, which we have not used before. For the purposes of this section, if $R$ is a cofibrant commutative $\S$-algebra, $x\in \pi_*R$ and $M$ is an $R$-module, then we define the $R[x^{-1}]$-module $M[x^{-1}]$ by $R[x^{-1}]\wedge_R M$ \cite[VII.4]{ekmm}.
\bprf Basterra \cite[Proposition 4.2]{basterra} proved that, if $A\to B\to C$ are maps of cofibrant commutative $\S$-algebras, then
\[TAQ(B|A) \wedge_B C \to TAQ(C|A) \to TAQ(C|B)\]
is a homotopy cofiber sequence of $C$-modules. Recall from (\ref{defku}) that we defined $KU$ as $Q\S[K(\Z,2)][x^{-1}]$, where $Q$ is a cofibrant replacement functor in $\SCAlg$ and $x\in \pi_2\S[K(\Z,2)]$. The following sequence of cofibrant commutative $\S$-algebras
\[\S\to Q\S[K(\Z,2)]\to Q\S[K(\Z,2)][x^{-1}]\]
begets a homotopy cofiber sequence of $KU$-modules
\begin{equation}\label{hocotaq}TAQ(Q\S[K(\Z,2)]) \wedge_{Q\S[K(\Z,2)]} KU \to TAQ(KU) \to TAQ(KU|Q\S[K(\Z,2)]).\end{equation}
Now, %
$TAQ(KU|Q\S[K(\Z,2)])$ is contractible, since $Q\S[K(\Z,2)]\to KU$ is a localization map \cite[Remark 3.4]{mccarthy-minasian}. Since by definition the leftmost factor of (\ref{hocotaq}) is $TAQ(Q\S[K(\Z,2)])[x^{-1}]$, 
the sequence (\ref{hocotaq}) gives a weak equivalence of $KU$-modules
\begin{equation}\label{turka}TAQ(Q\S[K(\Z,2)])[x^{-1}]\stackrel{\sim}{\to} TAQ(KU).\end{equation}
But \cite[Theorem 5]{basterra-mandell} gives that if $G$ is a topological abelian group, then the $Q\S[G]$-modules $TAQ(Q\S[G])$ and $Q\S[G] \wedge \mathbf{G}$ are weakly equivalent. 
Taking $G=K(\Z,2)$, localizing this equivalence at $x$ and combining it with (\ref{turka}), we get weak equivalences of $KU$-modules
\[KU \wedge \K \simeq TAQ(Q\S[K(\Z,2)])[x^{-1}] \stackrel{\sim}{\to} TAQ(KU). \qedhere\]
\eprf

We thank the anonymous referee for pointing us in the direction of the proof of the following fact, of which our previous proof was less simple.

\bcor \label{cortaq} The $KU$-modules $TAQ(KU)$ and $KU_\Q$ are weakly equivalent.
\bprf From Bott periodicity and the comment following the statement of Theorem \ref{taqku1}, we get that $TAQ(KU)$ and $KU\wedge H\Z$ are weakly equivalent. But the map $KU \wedge H\Z\to KU \wedge H\Q$ induced from the inclusion $\Z\subset \Q$ is a weak equivalence \cite[16.25]{switzer}, hence the result. %
\eprf
\ecor

\bobs \phantomsection \label{jurs} \be
\item \label{taq-od} The topological André-Quillen $B$-module $TAQ(B|A)$ can be computed as a stabilization, as follows from the work of \cite{basterra-mandell} and as made more explicit e.g. in \cite[Page 164]{schlichtkrull-higher}. More precisely, there is a tower with $\Omega^n (S^n\tilde\otimes_A B)$ in level $n$ whose homotopy colimit is weakly equivalent to $TAQ(B|A)$; here $S^n\tilde\otimes_A B$ is the $B$-module which is the cofiber of the map $B\to S^n\otimes_A B$ given by the inclusion of the basepoint in $S^n$. The symbol $\otimes_A$ denotes the tensor over $\Top$ of the category of commutative $A$-algebras. From Theorem \ref{holo} we deduce that $S^n\tilde\otimes KU \simeq KU \wedge K(\Z,n+2)$: we have identified these in Proposition \ref{lemgen}. This indicates a different way of computing $TAQ(KU)$.

\item \label{taq-ku-thom} Compare Theorem \ref{taqku1} with the reformulation %
found e.g. in \cite[Page 164]{schlichtkrull-higher} of a result of \cite{basterra-mandell}. %
It states that if $f:X\to BF$ is a map of $\infty$-loop spaces where $BF$ is a classifying space for stable spherical fibrations, then $TAQ(T(f))\simeq T(f)\wedge \mathbf{X}$. Just as in Remark \ref{x-thom}, the result for $TAQ(KU)$ coincides with the result we would obtain if we knew that $KU$ was somehow the Thom spectrum of a map $K(\Z,2)\to BU$.

\item \label{hokr} Consider the version of the Hochschild-Kostant-Rosenberg theorem in \cite[Theorem 1.1]{mccarthy-minasian}: if $A$ is a connective smooth commutative $\S$-algebra, there is a weak equivalence of commutative $A$-algebras $F(\Sigma TAQ(A))\stackrel{\sim}{\to}THH(A)$, where $F:A\Mod\to A\CAlg$ is the free commutative algebra functor.\footnote{Note, however, that the proof contains a gap: see \cite[Footnote 4]{antieau-vezzosi}. We thank Benjamin Antieau for pointing this out to us.} This statement does not apply to $KU$ since $KU$ is not connective (we have not checked the smoothness condition), but the conclusion is true (Theorem \ref{thhku2} and Corollary \ref{cortaq}). Just as in Remarks \ref{thh-thom}, \ref{x-thom} and the one just above, here is an example of a theorem that does not apply to $KU$ because $KU$ is not connective, but whose conclusion is nonetheless true. We speculate that there should be a version of the HKR theorem for $E_\infty$-ring spectra which dispenses with the connectiveness hypothesis. I would like to thank Tomasz Maszczyk for asking me about the HKR theorem in relation to Theorem \ref{thhku2}, thus inciting me to make these reflections. %
\ee
\eobs

\bibliographystyle{alpha}
\bibliography{../../../INCLUDE/main.bib}
\end{document}